\documentclass{amsart}
\usepackage{amscd}
\usepackage{amsthm}
\usepackage{amsmath}
\usepackage{graphicx}
\usepackage{amsfonts}
\usepackage{amssymb}
\theoremstyle{plain} 
\newtheorem{theorem}{Theorem}[section]
\newtheorem{corollary}[theorem]{Corollary}
\newtheorem{lemma}[theorem]{Lemma}
\newtheorem{proposition}[theorem]{Proposition}

\newtheorem{remark}{Remark}[section]
\newtheorem{ejemplo}{Example}[section]

\begin{document}
\title[The moduli of $G$-structures with linear connection.]{On the structure of the moduli of jets of $G$-structures with linear connection.}
\author{C. Mart\'{\i}nez Ontalba}
\address{Departamento de Geometr\'{\i}a y Topolog\'{\i}a\\
Escuela Superior de Inform\'{a}tica\\
Universidad Complutense de Madrid\\
28040-Madrid\\
Spain}
\email{celia\_martinez@mat.ucm.es}
\author{J. Mu\~{n}oz Masqu\'{e}}
\address{Instituto de F\'{\i}sica Aplicada\\
CSIC\\
Serrano 144\\
28006-Madrid\\
Spain}
\email{jaime@iec.csic.es}
\author{A. Vald\'{e}s}
\address{Departamento de Geometr\'{\i}a y Topolog\'{\i}a\\
Facultad de Ciencias Matem\'{a}ticas\\
Universidad Complutense de Madrid\\
28040-Madrid\\
Spain}
\email{avaldes@eucmos.sim.ucm.es}
\thanks{Supported by DGCYT under grant PB95-0124.}
\thanks{This paper is in final form and no version of it will be submitted for
publication elsewhere.}
\date{June, 1998}
\maketitle
\begin{abstract}The moduli space of jets of certain $G$-structures (basically those which
admit a canonical linear connection) is shown to be isomorphic to the quotient
of a natural $G$-module by $G$.
\end{abstract}

\section{Introduction}

Given a closed subgroup $G$ of the general linear group $GL(n,\mathbb{R)}$, a
$G$-structure is a reduced bundle $P(M,G)$ with structure group $G$ of the
bundle of linear frames $FM\rightarrow M$.

The main types of geometries arise from different choices of $G$. For example,
there is a one-to-one correspondence between the set of Riemannian metrics on
$M$ and the set of $O(n)$-structures on $M$. Analogously, almost Hermitian
geometries correspond to $U(n/2)$-structures, almost symplectic geometries to
$Sp(n/2)$-structures, and so on.

It is a well-known fact that $G$-structures on $M$ are in one-to-one
correspondence with smooth sections $s \in\Gamma(FM/G)$ of the quotient bundle
$FM/G \to M$. Denoting by $q:FM \to FM/G$ the natural projection, the fiber
over each $x \in M$ of the $G$-structure $P_{s} \to M$ associated to $s
\in\Gamma(FM/G)$ is $\left(  P_{s}\right)  _{x} = \{u=\left(  x;X_{1}, \ldots,
X_{n}\right)  \in F_{x}M \mid q(u) = s(x)\}$.

The group $\mathrm{Diff}(M)$ of diffeomorphisms of $M$ acts in a natural way
on the set of $G$-structures as follows. For each diffeomorphism
$f:M\rightarrow M^{\prime}$, there is an associated diffeomorphism $\bar
{f}:FM\rightarrow FM^{\prime}$ given by $\bar{f}\left(  x;X_{1},\ldots
,X_{n}\right)  =\left(  f(x);f_{\star}(X_{1}),\ldots,f_{\star}(X_{n})\right)
$, which defines a diffeomorphism $\tilde{f}:FM/G\rightarrow FM/G$ by
$\tilde{f}[u]=[\bar{f}(u)]$. The action of $\mathrm{Diff}(M)$ on
$\Gamma(FM/G)$ is defined by:
\[
f\cdot s=\tilde{f}\circ s\circ f^{-1},\hspace{0.5cm}f\in\mathrm{Diff}
(M),~s\in\Gamma(FM/G).
\]

Two $G$-structures $s$ and $s^{\prime}$ are said to be equivalent if they are
related by a diffeomorphism $f\in\mathrm{Diff}(M)$, which amounts to the fact
that $\bar{f}(P_{s})=P_{s^{\prime}}$, and they are said to be locally
equivalent at points $p\in M$ and $p^{\prime}\in M^{\prime}$ if they are
equivalent in some open neighborhoods of $p$ and $p^{\prime}$ by a
diffeomorphism which maps $p$ to $p^{\prime}$. We call the quotient
$\frak{M}_{G}(M)=\Gamma(FM/G)/\mathrm{Diff}(M)$ the \emph{moduli space of
$G$-structures} on $M$. The description of this space is a basic problem in
differential geometry.

It is possible to state an analogous problem for analytic $G$-structures. To
study this category, it is natural to introduce the spaces $J^{r}\left(
FM/G\right)  $ of jets of $G$-structures. The action of $\mathrm{Diff}(M)$ on
$\Gamma(FM/G)$ induces a natural action of the groupoid $J_{\mathrm{inv}%
}^{r+1}(M,M)$ of $(r+1)$-jets of diffeomorphisms of $M$ on the space
$J^{r}\left(  FM/G\right)  $ as follows:
\[
\left(  j_{x}^{r+1}f\right)  \cdot\left(  j_{x}^{r}s\right)  =j_{f(x)}%
^{r}\left(  f\cdot s\right)  ,\hspace{0.5cm}j_{x}^{r+1}f\in J_{\mathrm{inv}%
}^{r+1}(M,M),~j_{x}^{r}s\in J^{r}\left(  FM/G\right)  .
\]
We will call the quotient $\frak{M}_{G}^{r}(M)=J^{r}\left(  FM/G\right)
/J_{\mathrm{inv}}^{r+1}(M,M)$ the \emph{moduli space of $r$-jets of
$G$-structures} on $M$. There are natural projections $\frak{M}_{G}%
^{r+k}(M)\rightarrow\frak{M}_{G}^{r}(M)$, $k\geq0$, so that we can define the
moduli space of jets of $G$-structures as the projective limit $\frak{M}%
_{G}^{\infty}(M)=\underset{\leftarrow}{\lim}\,\frak{M}_{G}^{r}(M).$ The local
equivalence problem for analytic $G$-structures can be reduced to the study of
this moduli space.

Another reason for analyzing $\frak{M}_{G}^{\infty}(M)$ is that it is the
space where the geometric objects associated to $G$-structures are defined.

In general, the description of $\frak{M}_{G}^{\infty}(M)$ is an extremely
complicated problem. The aim of this paper is to describe, to some extent, the
structure of this space in the particular case of $G\subset GL(n,\mathbb{R)}$
being a closed subgroup such that the first prolongation $\frak{g}^{(1)}$
vanishes and there exists a supplementary $G$-submodule\footnote{Throughout
this paper, by $G$-module we mean a linear representation of $G$ in a finite
dimensional vector space.} $W$ of $\delta\left(  \left(  \mathbb{R}%
^{n}\right)  ^{\star}\otimes\frak{g}\right)  $ in $\bigwedge\nolimits^{2}%
\left(  \mathbb{R}^{n}\right)  ^{\star}\otimes\mathbb{R}^{n}.$ As we will
explain later, this is a technical condition which assures that $G$-structures
have a canonical linear connection attached.

More precisely, we will prove the following theorem.

\begin{theorem}
\label{Theorem}Let $M$ be a $n$-dimensional smooth manifold, and let $G\subset
GL\left(  n,\mathbb{R}\right)  $ be a closed subgroup such that the first
prolongation $\frak{g}^{(1)}$ of its Lie algebra $\frak{g}$ vanishes and there
exists a supplementary $G$-submodule $W$ of $\delta\left(  \left(
\mathbb{R}^{n}\right)  ^{\star}\otimes\frak{g}\right)  $ in $\bigwedge
\nolimits^{2}\left(  \mathbb{R}^{n}\right)  ^{\star}\otimes\mathbb{R}^{n}$.
Then, there exists a family of $G$-modules $\mathbf{S}^{r}$ and homomorphisms
$\mathbf{S}^{r+k}\rightarrow\mathbf{S}^{r}$, $k\geq0$, such that each space
$\frak{M}_{G}^{r}(M)$ of $r$-jets of $G$-structures is canonically isomorphic
to the quotient $\mathbf{S}^{r}/G$.

The moduli space $\frak{M}_{G}^{\infty}(M)$ is then canonically isomorphic to
the quotient $\mathbf{S}^{\infty}/G,$ where $\mathbf{S}^{\infty}%
=\underset{\leftarrow}{\lim}\,\mathbf{S}^{r}$.
\end{theorem}

\section{Natural linear connections}

In this section we will recall a classical result (\cite{B}) which shows that,
under the conditions on $G$ stated in Theorem \ref{Theorem}, it is possible to
associate a linear connection to each $G$-structure in a natural way. We will
also state certain properties of this assignment which will be useful later.

Let $G\subset GL(n,\mathbb{R})$ be a closed Lie subgroup, and $\frak{g}\subset
gl(n,\mathbb{R})$ its Lie algebra. Let us denote $V=\mathbb{R}^{n}$. Under the
identification $gl(n,\mathbb{R})\cong V^{\star}\otimes V$, we can consider
$\frak{g}$ as a Lie subalgebra of the endomorphisms of $V,$ $\mathrm{End}%
(V)\cong V^{\star}\otimes V$.

The alternation operator $\delta:V^{\star}\otimes\frak{g}\rightarrow
\bigwedge^{2}V^{\star}\otimes V$ is given by
\[
\delta\tau(u,v)=\tau(u)v-\tau(v)u,\hspace{0.5cm}\tau\in V^{\star}
\otimes\frak{g},u,v\in V.
\]
The first prolongation of $\frak{g}$ is then defined as the vector space
$\frak{g}^{(1)}:=\ker\delta$ and, from now on, we will assume that this
prolongation vanishes, i.e., $\frak{g}^{(1)}=\{0\}$.

The natural action of $G$ on $V$ defines linear $G$-actions on the spaces
$V^{\star}\otimes\frak{g}$ and $\bigwedge^{2}V^{\star}\otimes V$, with respect
to which the operator $\delta$ is an homomorphism of $G$-modules. The image of
this operator is a $G$-submodule and we assume that it admits a supplementary
$G$-submodule $W$:
\begin{equation}
\mbox{$\bigwedge\nolimits^{2}V^{\star}\otimes V$}=\delta\left(  V^{\star
}\otimes\frak{g}\right)  \oplus W. \label{supplementary}%
\end{equation}

Given a $G$-structure $s\in\Gamma(FM/G)$, the bundle $\bigwedge^{2}T^{\star
}M\otimes TM$ is identified with the associated bundle $P_{s}\times_{G}\left(
\bigwedge^{2}V^{\star}\otimes V\right)  $. The decomposition
(\ref{supplementary}) allows to express this bundle as the sum:
\[
\mbox{$\bigwedge\nolimits^{2}T^{\star}M\otimes TM$}=\left(  P_{s}\times_{G}
\delta\left(  V^{\star}\otimes\frak{g}\right)  \right)  \oplus\left(
P_{s}\times_{G}W\right)  .
\]

Let $\theta:T\left(  FM\right)  \rightarrow V$ be the canonical 1-form of
$FM$, which maps each tangent vector $X_{u}$ at a frame $u\in FM$ to
$\theta\left(  X_{u}\right)  =u^{-1}\left(  \pi_{\star}\left(  X_{u}\right)
\right)  \in V,$ where $\pi:FM\rightarrow M$ is the canonical projection and
the frame $u$ is understood as a linear isomorphism $u:V\rightarrow
T_{\pi\left(  u\right)  }M.$ We also denote by $\theta$ the restriction of the
canonical form to the subbundle $P_{s}.$ The restriction of the projection
$\pi_{\star}$ to each horizontal complement $H_{u}\subset T_{u}P_{s}$ to the
vertical tangent space at $u\in P_{s}$ is an isomorphism $H_{u}\cong
T_{\pi\left(  u\right)  }M$, so that for any frame $u$ and any vector $v\in V$
there is a unique $B_{u}\left(  v\right)  \in H_{u}$ such that
%
%
%
%
$\theta\left(  B_{u}\left(  v\right)  \right)  =v.$ The horizontal complement
$H_{u}$ defines then an element $t(H_{u})\in\bigwedge^{2}V^{\star}\otimes V$
by:
\[
t(H_{u})(v,w)=d\theta\left(  B_{u}(v),B_{u}(w)\right)  ,~~~v,w\in V.
\]

If $H_{u}$ are the horizontal subspaces of a linear connection on $M$ adapted
to $P_{s}$, then there is a well-defined and $G$-equivariant \emph{torsion
function}: $P_{s}\rightarrow\bigwedge^{2}V^{\star}\otimes V$, given by
$u\mapsto t(H_{u})$, which is related to the torsion tensor $T\in\Gamma\left(
\bigwedge^{2}T^{\star}M\otimes TM\right)  $ by: $t(H_{u})\left(  v,w\right)
=u^{-1}\left(  T_{x}\left(  u\left(  v\right)  ,u\left(  w\right)  \right)
\right)  $ for each $u\in P_{s}$ and $v,w\in V.$

It can be shown (see, e.g., \cite{B} and \cite{F}) that the condition
$t(H_{u}) \in W$ characterizes the horizontal subspaces defining a linear
connection adapted to $P_{s}$. In other words:

\begin{theorem}
[\cite{B}]\label{Bernard}Let $G\subset GL(n,\mathbb{R})$ be a closed Lie
subgroup
%
%
%
%
such that the first prolongation $\frak{g}^{(1)}$ vanishes and there exists a
supplementary $G$-submodule $W$ of $\delta\left(  V^{\star}\otimes
\frak{g}\right)  $ in $\bigwedge\nolimits^{2}V^{\star}\otimes V.$ Then, for
each $G$-structure $s\in\Gamma(FM/G)$ and each $u\in P_{s}:=s\left(  M\right)
$ there exists a unique horizontal space $H_{u}\subset T_{u}P$ such that
$t(H_{u})\in W.$ These horizontal spaces define a linear connection
$\nabla(s)$ adapted to $P_{s}$ which is characterized by the condition that
its torsion is a section of the vector bundle $P_{s}\times_{G}W$.
\end{theorem}

Next two lemmas provide some properties of this \emph{canonical connection}
which will be used to describe the structure of the moduli spaces in later sections.

\begin{lemma}
\label{order1} The assignment $s\rightarrow\nabla(s)$ defines an operator
$\nabla:\Gamma(FM/G)\rightarrow\Gamma(C(M))$, with $C(M)\rightarrow M$ being
the affine bundle of linear connections on $M$, which satisfies:

\begin{enumerate}
\item $\nabla$ is natural, i.e., for each diffeomorphism $f:M\rightarrow
M^{\prime}$, the direct image of $\nabla(s)$ onto $\bar{f}(P_{s})$ is
$\nabla(f\cdot s)$.

\item $\nabla(s)$ is an operator of order 1, i.e., if $j_{x}^{1}s=j_{x}%
^{1}s^{\prime}$, then $\nabla(s)(x)=\nabla(s^{\prime})(x)$. Thus, $\nabla$
defines a map, which we will also denote $\nabla:J^{1}(FM/G)\rightarrow C(M)$,
by $\nabla(j_{x}^{1}s)=\nabla(s)(x)$.
\end{enumerate}
\end{lemma}

\textbf{Proof.} In order to prove the naturality of the map $\nabla
:\Gamma(FM/G)\rightarrow\Gamma(C(M))$, it suffices to show that $t(H_{u})\in
W$ implies $t\left(  \bar{f}_{\star}(H_{u})\right)  \in W$ for each
diffeomorphism $f:M\rightarrow M^{\prime}$ because, then, the direct image
$\bar{f}_{\star}(H_{u})\subset\bar{f}_{\star}\left(  T_{u}P\right)
=T_{\bar{f}\left(  u\right)  }P_{f\cdot s}$
%
%
%
%
of each horizontal space of the connection $\nabla(s)$ by $f$ is, by Theorem
\ref{Bernard}, a horizontal space of the connection $\nabla\left(  f\cdot
s\right)  $.
%
%
%
%
But this follows easily from the functoriality of the canonical 1-form of the
bundle of linear frames. Explicitly, the canonical 1-form $\theta^{\prime
}:TFM^{\prime}\rightarrow V$ of $FM^{\prime}$ is related to $\theta$ by
$\theta=\bar{f}^{\star}\theta^{\prime}$ as one can easily check. From this it
follows that $\theta^{\prime}\left(  \bar{f}_{\star}(B_{u})(v)\right)
=\theta\left(  (B_{u})(v)\right)  =v,$ for all $v\in V,$ and therefore we have
that $B_{\bar{f}\left(  u\right)  }^{\prime}=f_{\star}\left(  B_{u}\right)  $
and so%

\begin{align*}
t\left(  \bar{f}_{\star}(H_{u})\right)  (v,w)  &  =d\theta^{\prime}\left(
B_{\bar{f}\left(  u\right)  }^{\prime}(v),B_{\bar{f}\left(  u\right)
}^{\prime}(w)\right) \\
&  =d\theta^{\prime}\left(  \bar{f}_{\star}(B_{u})(v),\bar{f}_{\star}%
(B_{u})(w)\right) \\
&  =d\left(  \bar{f}^{\star}\theta^{\prime}\right)  \left(  (B_{u}%
)(v),(B_{u})(w)\right) \\
&  =d\theta\left(  (B_{u})(v),(B_{u})(w)\right)  =t(H_{u})(v,w).
\end{align*}

To prove that $\nabla$ is an operator of order 1, first notice that the
condition $j_{x}^{1}s=j_{x}^{1}s^{\prime}$ is equivalent to the fact that the
mappings $s_{\star}:T_{x}M\longrightarrow T_{s\left(  x\right)  }\left(
FM/G\right)  $ and $s_{\star}^{\prime}:T_{x}M\longrightarrow T_{s^{\prime
}\left(  x\right)  }\left(  FM/G\right)  $ coincide. But, then, the equality
$T_{u}P_{s} =T_{u}P_{s^{\prime}}$ holds for any $u\in s(x)=s^{\prime}(x)$. To
see this, notice that if $q:FM\rightarrow FM/G$ is, as before, the natural
projection, then:%

\begin{align*}
T_{u}P_{s}  &  =q_{\star}^{-1}\left(  T_{q(u)}\left(  s(M)\right)  \right) \\
&  =q_{\star}^{-1}\left(  s_{\star}\left(  T_{x}M\right)  \right) \\
&  =q_{\star}^{-1}\left(  s_{\star}^{\prime}\left(  T_{x}M\right)  \right) \\
&  =T_{u}P_{s^{\prime}}.
\end{align*}
Since $t(H_{u})\in W$ is an algebraic condition depending only on the vector
spaces $T_{u}P_{s}$ and $W$, it follows that the equality $T_{u}P_{s}
=T_{u}P_{s^{\prime}}$ implies the equality of the horizontal subspaces
defining the connections: $H_{u}=H_{u}^{\prime}$ for $u\in s(x)=s^{\prime}
(x)$.$\blacksquare$

Now, let $\sigma:\mathcal{U}\subset M \to P_{s}$ be a section defined in a
coordinate neighborhood $\left(  \mathcal{U},(x^{1}, \ldots, x^{n})\right)  $
of $x\in M$. We will denote $\sigma= \left(  X_{1},\ldots, X_{n}\right)  $,
with $X_{i} = \sum_{i}\sigma_{ji}\frac{\partial}{\partial x^{j}}$ for each $i
= 1, \ldots n$. Let $\omega\in\Gamma\left(  T^{\star}P_{s}\otimes
\frak{g}\right)  $ be the connection form of $\nabla:=\nabla(s)$.

If we consider, as before, each $\sigma(x) \in P_{s}$ as a linear isomorphism
$\sigma(x):V\to T_{x}M$, we can define a function $\eta\in C^{\infty}\left(
\mathcal{U}, V^{\star}\otimes\frak{g}\right)  $ by: $\eta(x) = (\sigma^{\star
}\omega)_{x}\circ\sigma(x)$. The components of $\eta$ in the standard basis
$\{v_{1}, \ldots, v_{n}\}$ of $V$ are related to those of the local connection
1-form $\sigma^{\star}\omega\in\Gamma\left(  T^{\star}M\otimes\frak{g}\right)
$ by: $\eta^{k}_{ij} = \left(  \sigma^{\star}\omega\right)  ^{k}_{j}(X_{i})$.

Let $P_{\mathrm{im}\delta}:\bigwedge^{2}V^{\star}\otimes V\rightarrow
\delta\left(  V^{\star}\otimes\frak{g}\right)  $ be the projection onto
$\delta\left(  V^{\star}\otimes\frak{g}\right)  $ according to the
decomposition (\ref{supplementary}). Since $\ker\delta=\{0\}$, it makes sense
to consider the homomorphism:
\[
\delta^{-1}\circ P_{\mathrm{im}\delta}:\bigwedge\nolimits^{2}V^{\star}\otimes
V\rightarrow V^{\star}\otimes\frak{g}.
\]

\begin{lemma}
The local connection 1-form $\sigma^{\star}\omega\in\Gamma\left(  T^{\star
}M\otimes\frak{g}\right)  $ is determined by the equation: $\eta(x)=-\left(
\delta^{-1}\circ P_{\mathrm{im}\delta}\right)  \left(  \tilde{t}%
(\sigma(x))\right)  ,$ where $\tilde{t}$ is the torsion function of the linear
flat connection making parallel $\sigma$.
\end{lemma}

\textbf{Proof. }Using the identity $\nabla_{X}X_{j}=\sum_{k}X_{k}\left(
\sigma^{\star}\omega\right)  _{j}^{k}(X),$ $X\in\Gamma(TM),~1\leq j\leq n,$ we
see that the components of the torsion tensor of $\nabla$ in the moving frame
$\sigma,$ which are given by $T(X_{i},X_{j})=\nabla_{X_{i}}X_{j}-\nabla
_{X_{j}}X_{i}-[X_{i},X_{j}]=\sum_{k}T_{ij}^{k}X_{k}$, $1\leq i,j\leq n,$ can
be written as $T_{ij}^{k}(x)=\eta_{ij}^{k}(x)-\eta_{ji} ^{k}(x)+\tilde{t}%
_{ij}^{k}(\sigma(x)).$ Using the relation $t(H_{\sigma(x)})(v_{i},v_{j}%
)=\sum_{k}T_{ij}^{k}(x)v_{k}$, last equation reads:
\[
t(H_{\sigma(x)})=\delta\eta(x)+\tilde{t}(\sigma(x)).
\]
The condition $t(H_{\sigma(x)})\in W$ can be written as $\left(  \delta
^{-1}\circ P_{\mathrm{im}\delta}\right)  \left(  t(H_{\sigma(x)})\right)  =0$
or, equivalently, $\eta(x)=-\left(  \delta^{-1}\circ P_{\mathrm{im} \delta
}\right)  \left(  \tilde{t}(\sigma(x))\right)  $.$\blacksquare$

\begin{remark}
\label{forF}\emph{Notice that $\tilde{t}_{ij}^{k}(\sigma(x))=-\sum_{h,l}%
\sigma^{kl}\left(  \sigma_{hi}\frac{\partial\sigma_{lj}}{\partial x^{h}%
}-\sigma_{hj}\frac{\partial\sigma_{li}}{\partial x^{h}}\right)  $, where
$\left(  \sigma^{ij}\right)  $ stands for the inverse matrix of $\left(
\sigma_{ij}\right)  $. Thus, the components $\left(  \sigma^{\star}%
\omega\right)  _{\beta}^{\gamma}\in\Gamma(T^{\star}M)$ of the local connection
1-form $\sigma^{\star}\omega\in\Gamma\left(  T^{\star}M\otimes\frak{g}\right)
$ are given by:
\begin{equation}
\left(  \sigma^{\star}\omega\right)  _{\beta}^{\gamma}(X_{\alpha}%
)=\sum_{i,j,k,h,l}A_{\alpha\beta k}^{\gamma ij}\sigma^{kl}\left(  \sigma
_{hi}\frac{\partial\sigma_{lj}}{\partial x^{h}}-\sigma_{hj}\frac
{\partial\sigma_{li}}{\partial x^{h}}\right)  ,~~~1\leq\alpha,\beta,\gamma\leq
n,\label{local_connection}%
\end{equation}
where the real coefficients $A_{\alpha\beta k}^{\gamma ij}$ are determined by
the homomorphism $\delta^{-1}\circ P_{\mathrm{im}\delta}$. These coefficients
are then universal, in the sense that they only depend on the group $G\subset
GL(n,\mathbb{R})$ and the supplementary $W$ chosen to define the connections,
but not on the particular $G$-structure considered.}
\end{remark}

\begin{remark}
\label{smoothness} \emph{Equation (\ref{local_connection}) shows that the map
$\widetilde{\nabla}:J^{1}(FM)\rightarrow C(M)$ defined as $\widetilde{\nabla
}\left(  j_{x}^{1}\sigma\right)  =\nabla\left(  q\circ\sigma\right)  (x)$ is
smooth. Since the diagram
\[%
\begin{array}
[c]{rcl}%
J^{1}(FM) & \overset{q^{1}}{\longrightarrow} & J^{1}(FM/G)\\
{\scriptstyle\widetilde{\nabla}}\searrow &  & \swarrow{\scriptstyle\nabla}\\
&  C(M) &
\end{array}
\]
is commutative, and the projection $q^{1}:J^{1}(FM)\rightarrow J^{1}(FM/G)$ is
a surjective submersion, we conclude that the map $\nabla:J^{1}%
(FM/G)\rightarrow C(M)$ is also smooth.}
\end{remark}

To end this section, we illustrate the construction of canonical connections
with some examples. In all cases we will take $W$ as the $G$-submodule of
tensors $T\in\bigwedge^{2}V^{\star}\otimes V$ such that trace $(A\circ
i_{v}T)=0$ for every $A\in\frak{g}$ and every $v\in V$. This is a
supplementary $G$-submodule of the image of $\delta$ whenever $\frak{g}%
^{(1)}=0$ and $\frak{g}$ is closed under transposition (see \cite{V}).

\begin{ejemplo}
\label{Levi-Civita} \emph{It is well-known that the canonical connection (in
the sense of Theorem \ref{Bernard}) associated to each $O(n)$-structure or,
equivalently, to each Riemannian metric is the Levi-Civita connection. That is
because for $G=O(n)$ the alternation operator $\delta$ is surjective and,
therefore, the supplementary submodule in equation (\ref{supplementary}) is
$W=\{0\}$. The condition $\eta(x)=-\delta^{-1}\left(  \tilde{t}(\sigma
(x))\right)  $ is equivalent to Koszul formula, as we next explain. }

\emph{The operator $\delta^{-1}\circ P_{\mathrm{im}\delta}=\delta
^{-1}:\bigwedge\nolimits^{2}V^{\star}\otimes V\rightarrow V^{\star}%
\otimes\frak{g}$ is easily obtained as follows. Let us denote by $\langle
\cdot,\cdot\rangle$ the standard inner product in $V$, and take any
$T=\delta\tau\in\bigwedge\nolimits^{2}V^{\star}\otimes V$, $\tau\in V^{\star
}\otimes\frak{g}$. Then, for any $u$, $v$, $w\in V$, the following relations
hold:
\begin{align*}
\langle T(u,v),w\rangle &  =\langle\tau(u)v,w\rangle-\langle\tau
(v)u,w\rangle\\
\langle T(w,u),v\rangle &  =\langle\tau(w)u,v\rangle-\langle\tau
(u)w,v\rangle\\
\langle T(w,v),u\rangle &  =\langle\tau(w)v,u\rangle-\langle\tau(v)w,u\rangle
\end{align*}
Adding the three equations, taking into account that $\frak{g}=\{A\in
V^{\star}\otimes V\mid\langle Au,v\rangle+\langle Av,u\rangle=0~~~\forall
u,v\in V\}$, we obtain the expression of $\tau=\delta^{-1}T$:
\begin{equation}
\langle\tau(u)v,w\rangle=\frac{1}{2}\left(  \langle T(u,v),w\rangle+\langle
T(w,u),v\rangle+\langle T(w,v),u\rangle\right)  ,\label{inv_delta}%
\end{equation}
for all $u,v,w\in V$. Hence, the condition $\eta(x)=-\delta^{-1}\left(
\tilde{t}(\sigma(x))\right)  $ can be written, in components, as:
\begin{equation}
\eta_{ij}^{k}(x)=-\frac{1}{2}\left(  \tilde{t}_{ij}^{k}(\sigma(x))+\tilde
{t}_{ki}^{j}(\sigma(x))+\tilde{t}_{kj}^{i}(\sigma(x))\right)  .\label{eta}%
\end{equation}
}

\emph{Now, let $(M,g)$ be a Riemannian manifold and $P_{s}\rightarrow M$ the
corresponding $O(n)$-structure, that is, the subbundle of orthonormal frames
with respect to $g$. For any section $\sigma=\left(  X_{1},\ldots
,X_{n}\right)  :\mathcal{U}\rightarrow P_{s}$, the components of $\eta\in
C^{\infty}\left(  \mathcal{U},V^{\star}\otimes\frak{g}\right)  $ are given by:
$\eta_{ij}^{k}=\left(  \sigma^{\star}\omega\right)  _{j}^{k}(X_{i}%
)=g(\nabla_{X_{i}}X_{j},X_{k})$, whereas those of $\tilde{t}(\sigma(x))$ are:
$\tilde{t}_{ij}^{k}=-g([X_{i},X_{j}],X_{k})$. Substitution of these
expressions in (\ref{eta}) yields:
\[
2g(\nabla_{X_{i}}X_{j},X_{k})=g([X_{i},X_{j}],X_{k})+g([X_{k},X_{i}%
],X_{j})+g([X_{k},X_{j}],X_{i}),
\]
which is nothing but Koszul formula in the orthonormal frame $\{X_{1},\ldots
X_{n}\}$. }
\end{ejemplo}

\begin{ejemplo}
\label{foliations}\emph{Consider now the case of $\left(  O(p)\times
O(q)\right)  $-structures, i.e., those with structure group $G=\{\left(
\begin{array}
[c]{cc}%
A_{1} & 0\\
0 & A_{2}%
\end{array}
\right)  \in O(n)\mid A_{1}\in O(p),A_{2}\in O(q)\}$. These structures
correspond to Riemannian almost product structures. }

\emph{Let $\{v_{1},\ldots,v_{n}\}$ be the canonical basis of $V$, and let us
denote $I_{1}=\{1,\ldots,p\}$, $I_{2}=\{p+1,\ldots,p+q\}$, $V_{\alpha
}=\mathrm{Span}\{v_{i}\}_{i\in I_{\alpha}}$, $\alpha=1,2$. It can be easily
seen that $W$ is given as the $G$-submodule of tensors $\tilde{T}\in
\bigwedge^{2}V^{\star}\otimes V$ such that:
\begin{equation}
\langle\tilde{T}(u,v),w\rangle=\langle\tilde{T}(u,w),v\rangle~~~~~\mathrm{if}%
~~~~~v,w\in V_{\alpha},~\alpha=1,2,\label{W_foliations}%
\end{equation}
for all $u\in V$. }

\emph{ }

\emph{In order to determine an expression of the operator $\delta^{-1}\circ
P_{\mathrm{im}\delta}$, let us decompose any $T\in\bigwedge^{2}V^{\star
}\otimes V$ as $T=\tilde{T}+\delta\tau$, with $\tilde{T}\in W$ and $\tau\in
V^{\star}\otimes\frak{g}$.}

\emph{If $v\in V_{\alpha}$, $w\in V_{\beta}$, with $\alpha\neq\beta$, then
$\langle\tau(u)v,w\rangle=0$ for each $u\in V$, because $\tau(u)$ belongs to
$\frak{g}=o(p)\oplus o(q)$.}

\emph{If $u,v,w\in V_{\alpha}$, $\alpha=1,2$, then $\langle\tilde
{T}(u,v),w\rangle=0$, as follows easily from (\ref{W_foliations}). In this
case, one can compute, as in the previous example:
\[
\langle\tau(u)v,w\rangle=\frac{1}{2}\left(  \langle T(u,v),w\rangle+\langle
T(w,u),v\rangle+\langle T(w,v),u\rangle\right)  .
\]
}

\emph{Finally, if $u\in V_{\alpha}$ and $v,w\in V_{\beta}$, $\alpha\neq\beta$,
we have:
\begin{align*}
\langle T(u,v),w\rangle &  =\langle\tau(u)v,w\rangle+\langle\tilde
{T}(u,v),w\rangle\\
\langle T(u,w),v\rangle &  =\langle\tau(u)w,v\rangle+\langle\tilde
{T}(u,w),v\rangle.
\end{align*}
Substracting the second equation from the first one, we obtain:
\[
\langle\tau(u)v,w\rangle=\frac{1}{2}\left(  \langle T(u,v),w\rangle-\langle
T(u,w),u\rangle\right)  .
\]
}

\emph{Thus, in components, we have:
\[
\left(  (\delta^{-1}\circ P_{\mathrm{im}\delta})(T)\right)  _{ij}%
^{k}=\left\{
\begin{array}
[c]{lll}%
\frac{1}{2}\left(  T_{ij}^{k}+T_{ki}^{j}+T_{kj}^{i}\right)   & \text{if} &
i,j,k\in I_{\alpha},~~~~~\alpha=1,2,\\
\frac{1}{2}\left(  T_{ij}^{k}-T_{ik}^{j}\right)   & \text{if} & i\in
I_{\alpha},~~~j,k\in I_{\beta},~~~~~\alpha\neq\beta,\\
0 & \text{if} & j\in I_{\alpha},~~~k\in I_{\beta},~~~~~\alpha\neq\beta.
\end{array}
\right.
\]
}

\emph{Let us consider now a Riemannian almost product manifold, i.e., a
Riemannian manifold $(M,g)$ and a tensor field $\phi\in\frak{T}_{1}^{1}(M)$
such that $\phi^{2}=\mathrm{id}$ and $g(\phi X,\phi Y)=g(X,Y)$ for any vector
fields $X,Y$ on $M$. The tensor $\phi$ gives rise to two mutually orthogonal
distributions $V_{1}$ and $V_{2}$, corresponding to its eigenvalues 1 and -1,
called the vertical and horizontal distributions, respectively. Denoting
$p=\mathrm{dim}V_{1}$, $q=\mathrm{dim}V_{2}$, the corresponding $\left(
O(p)\times O(q)\right)  $-structure is the subbundle of orthonormal frames
$(X_{1},\ldots,X_{n})$ such that $X_{1},\ldots,X_{p}$ are vertical vectors and
$X_{p+1},\ldots,X_{n}$ are horizontal vectors. }

\emph{From the expression of $\delta^{-1}\circ P_{\mathrm{im}\delta}$, we
obtain, as in the preceding example:
\[
g(\nabla_{X_{i}}X_{j},X_{k})=\frac{1}{2}\left(  g([X_{i},X_{j}],X_{k}%
)+g([X_{k},X_{i}],X_{j})+g([X_{k},X_{j}],X_{i})\right)
\]
if $i,j,k\in I_{\alpha}$, $\alpha=1,2$,
\[
g(\nabla_{X_{i}}X_{j},X_{k})=\frac{1}{2}\left(  g([X_{i},X_{j}],X_{k}%
)-g([X_{i},X_{k}],X_{j})\right)
\]
if $i\in I_{\alpha}$, $j,k\in I_{\beta}$, $\alpha\neq\beta$, and
\[
g(\nabla_{X_{i}}X_{j},X_{k})=0
\]
if $j\in I_{\alpha}$, $k\in I_{\beta}$, $\alpha\neq\beta$. }

\emph{The third equation implies that $\nabla_{X}$ leaves the vertical and
horizontal distributions invariant for all $X$. The first one gives:
\begin{align*}
g(\nabla_{X}Y,Z) &  =\frac{1}{2}\left(  X(g(Y,Z))+Y(g(Z,X))-Z(g(X,Y))\right.
\\
&  \left.  +g([X,Y],Z)+g([Z,X],Y)-g([Y,Z],X)\right)  ,
\end{align*}
if $X,Y,Z\in V_{\alpha}$, $\alpha=1,2$. Denoting by $\nabla^{\mathrm{LC}}$ the
Levi-Civita connection, this means that, if $X,Y\in V_{\alpha}$, then
$\nabla_{X}Y$ is the orthogonal projection of $\nabla_{X}^{\mathrm{LC}}Y$ onto
$V_{\alpha}$ according to the decomposition $TM=V_{1}\oplus V_{2}$. Finally,
the second equation yields:
\[
g(\nabla_{X}Y,Z)=\frac{1}{2}\left(  X(g(Y,Z))+g([X,Y],Z)-g([X,Z],Y)\right)  ,
\]
if $X\in V_{\alpha}$, $Y,Z\in V_{\beta}$, $\alpha\neq\beta$.}
\end{ejemplo}

\begin{ejemplo}
\label{webs} \emph{Next, we consider the case $G=\mathbb{R}^{\star}%
\cdot\mathrm{I_{n}}$ with $n\geq2$. If $n=2$ we obtain classical webs over
surfaces (see \cite{C}). }

\emph{The Lie algebra $\frak{g}$ of $G$ is generated by the identity matrix
$\mathrm{I_{n}}$, hence we can take:
\[
W=\{T\in\hbox{$\bigwedge^2$}V^{\star}\otimes V\mid\mathrm{trace}%
~(i_{v}T)=0~~~\forall v\in V\}.
\]
A short computation shows that the operator $\delta^{-1}\circ P_{\mathrm{im}%
\delta}$ is defined by:
\[
\left(  (\delta^{-1}\circ P_{\mathrm{im}\delta})(T)\right)  (v)=\left(
\frac{1}{n-1}\mathrm{trace}~(i_{v}T)\right)  I_{n},
\]
for all $T\in\bigwedge^{2}V^{\star}\otimes V$, $v\in V$. The condition
$\eta(x)=-\delta^{-1}\left(  \tilde{t}(\sigma(x))\right)  $ is written, in
components:
\[
\eta_{ij}^{k}(x)=-\frac{\delta_{j}^{k}}{n-1}\sum_{l}\tilde{t}_{il}^{l}%
(\sigma(x)).
\]
}

\emph{From this expression, it follows that the canonical connection is given,
in any local section $\sigma=(X_{1},\ldots,X_{n})$, by:
\[
\nabla_{X_{i}}X_{j}=\left(  \frac{1}{n-1}\sum_{k}C_{ik}^{k}\right)  X_{j},
\]
with the coefficients $C_{ij}^{k}$ being defined by $[X_{i},X_{j}]=\sum
_{k}C_{ij}^{k}X_{k}$. }
\end{ejemplo}

\section{Canonical representation of the moduli spaces}

In this section we are going to describe the underlying manifolds of the
$G$-modules $\mathbf{S}^{r}$ to which we refer in Theorem \ref{Theorem}, and
the action of the group $G$ on them. However, the definition of linear
structures on these $G$-manifolds will be postponed until section 5. We also
describe the isomorphisms $\frak{M}_{G}^{r} (M)\cong\mathbf{S}^{r}/G$.

Let us denote by $\frak{G}_{0}^{r+1}\subset J_{\mathrm{inv}}^{r+1}\left(
V,V\right)  $ the Lie group of $\left(  r+1\right)  $-jets of diffeomorphisms
of $V$ which leave the origin $0\in V$ fixed. The restriction to $\frak{G}%
_{0}^{r+1}$ of the action of the Lie groupoid $J_{\mathrm{inv}}^{r+1}\left(
V,V\right)  $ on $J^{r}\left(  FV/G\right)  $  defines an action of
$\frak{G}_{0}^{r+1}$ on $J_{0}^{r}\left(  FV/G\right)  $.

\begin{proposition}
There exists a canonical bijection $\frak{M}_{G}^{r}(M)\cong J_{0}^{r}\left(
FV/G\right)  /\frak{G}_{0}^{r+1}.$
\end{proposition}

\textbf{Proof.} Let us recall that $\frak{M}_{G}^{r}(M)$ has been defined as
the quotient
\[
J^{r}\left(  FM/G\right)  /J_{\mathrm{inv}}^{r+1}\left(  M,M\right)  .
\]
For each element $[j_{x}^{r}s]\in\frak{M}_{G}^{r}(M)$, we choose a chart
$\varphi:\mathcal{U}\subset M\rightarrow V$, centered at $x\in M$, and define
the element $\left(  j_{x}^{r+1}\varphi\right)  \cdot\left(  j_{x}%
^{r}s\right)  =j_{0}^{r}\left(  \tilde{\varphi}\circ s\circ\varphi
^{-1}\right)  \in J_{0}^{r}\left(  FV/G\right)  $. If we take another
representative $\left(  j_{x}^{r+1}f\right)  \cdot\left(  j_{x}^{r}s\right)  $
of $[j_{x}^{r}s]$, and a chart $(\mathcal{U}^{\prime},\varphi^{\prime})$
centered at $f(x)\in M$, we have:
\[
\left(  j_{f(x)}^{r+1}\varphi^{\prime}\right)  \cdot\left(  \left(
j_{x}^{r+1}f\right)  \cdot\left(  j_{x}^{r}s\right)  \right)  =\left(
j_{0}^{r+1}(\varphi^{\prime}\circ f\circ\varphi^{-1})\right)  \cdot\left(
\left(  j_{x}^{r+1}\varphi\right)  \cdot\left(  j_{x}^{r}s\right)  \right)  .
\]
Thus, the element of $J_{0}^{r}\left(  FV/G\right)  $ assigned to $[j_{x}%
^{r}s]$ is determined up to an $(r+1)$-jet of the form $j_{0}^{r+1}%
(\varphi^{\prime}\circ f\circ\varphi^{-1})$, which is an element of
$\frak{G}_{0}^{r+1}$ and, hence, there is a well defined mapping $\frak{M}%
_{G}^{r}(M)\longrightarrow J_{0}^{r}\left(  FV/G\right)  /\frak{G}_{0}^{r+1}.$
It is easy to find the inverse of this mapping. Given a class $[j_{0}^{r}t]\in
J_{0}^{r}\left(  FV/G\right)  /\frak{G}_{0}^{r+1}$ and a chart $\varphi
:\mathcal{U}\subset M\rightarrow V$ centered at $x\in M$ the element $\left[
\left(  j_{0}^{r+1}\varphi^{-1}\right)  \cdot\left(  j_{0}^{r}t\right)
\right]  \in J^{r}\left(  FM/G\right)  /J_{\mathrm{inv}}^{r+1}\left(
M,M\right)  $ is well-defined. The mapping $[j_{0}^{r}t]\mapsto\left[  \left(
j_{0}^{r+1}\varphi^{-1}\right)  \cdot\left(  j_{0}^{r}t\right)  \right]  $ is
the required inverse $J_{0}^{r}\left(  FV/G\right)  /\frak{G}_{0}%
^{r+1}\longrightarrow\frak{M}_{G}^{r}(M).\blacksquare$

From now on, we will denote the spaces $\frak{M}_{G}^{r}(M)$ just by
$\frak{M}_{G}^{r}$, because we have seen that they actually do not depend on
the base manifold $M$.

Now, we define the spaces $\mathcal{E}^{r}(V)$ of \emph{framed $r$-jets of
$G$-structures} as:
\[
\mathcal{E}^{r}(V) := J_{0}^{r}\left(  FV/G\right)  \times_{F_{0}V/G} F_{0}V=
\{(j_{0}^{r}s,u)\in J^{r}_{0}\left(  FV/G\right)  \times F_{0}V:[u]=s\left(
0\right)  \}.
\]

We consider the left action of $\frak{G}_{0}^{r+1}$ on $\mathcal{E}^{r}(V)$
given by:
\[
\left(  j_{0}^{r+1}f\right)  \cdot\left(  j_{0}^{r}s,u\right)  =\left(
\left(  j_{0}^{r+1}f\right)  \cdot\left(  j_{0}^{r} s\right)  ,\bar{f}\left(
u\right)  \right)  ,\hspace{.5cm}j_{0}^{r+1}f \in\frak{G}_{0}^{r+1},~ \left(
j_{0} ^{r}s,u\right)  \in\mathcal{E}^{r}(V).
\]

Besides, we consider the right action of $G$ on $\mathcal{E}^{r}(V)$ induced
by the action of $G$ on $F_{0}V$, i.e.,
\[
\left(  j_{0}^{r}s,u\right)  \cdot g=\left(  j_{0}^{r}s,u\cdot g\right)
,\hspace{.5cm}g\in G,~\left(  j_{0}^{r}s,u\right)  \in\mathcal{E}^{r}(V).
\]
The quotient $\mathcal{E}^{r}(V) /G$ can obviously be identified with
$J_{0}^{r}\left(  FV/G\right)  $.

It is easy to check that the group actions just defined commute. Therefore, we
have the bijections:
\[
\frak{M}^{r}_{G} \cong\frac{J_{0}^{r}\left(  FV/G\right)  }{\frak{G}_{0}^{r+1}
}\cong\frac{\mathcal{E}^{r}(V)/G}{\frak{G}_{0}^{r+1}} \cong\frac
{\mathcal{E}^{r}(V)/\frak{G}_{0}^{r+1}}{G}.
\]

Thus, if there is a manifold $\mathbf{S}^{r}$ such that $\mathbf{S}
^{r}/G\cong\frak{M}_{G}^{r},$ the natural candidate is $\mathbf{S}
^{r}:=\mathcal{E}^{r}(V)/\frak{G}_{0}^{r+1}$, if it is indeed a smooth
manifold without singularities. In general, this would not be true. However,
the existence of linear connections functorially attached to $G$-structures
makes each $\mathcal{E}^{r}(V)$ a trivial principal bundle with structure
group $\frak{G}_{0}^{r+1}$, as we will see below, and so $\mathcal{E}%
^{r}(V)/\frak{G}_{0}^{r+1}$ is a manifold in our particular cases.

Given $(j_{0}^{r}s,u)\in\mathcal{E}^{r}\left(  V\right)  $, let $s:\mathcal{U}
\subset V\rightarrow FV/G$ be a representative of $j_{0}^{r}s$, defined on an
open neighborhood $\mathcal{U}$ of $0\in V$. The connection $\nabla(s)$
provides an exponential mapping $\exp_{s}:\mathcal{W}_{0}\subset
T_{0}V\rightarrow V$ defined on some neighborhood $\mathcal{W}_{0}$ of $0\in
T_{0}V$. The composition of this mapping with the isomorphism $u:V\to T_{0}V$
yields a diffeomorphism defined on a neighborhood of $0\in V$, which is
nothing but the set of normal coordinates associated to the connection
$\nabla(s)$ and the frame $u$.

\begin{lemma}
The assignment $\mathrm{Exp}^{r}:\mathcal{E}^{r}\left(  V\right)
\rightarrow\frak{G}_{0}^{r+1}$ given by $\mathrm{Exp}^{r}\left(  j_{0}%
^{r}s,u\right)  =j_{0}^{r+1}\left(  \exp_{s}\circ u\right)  $ is a
well-defined smooth map. Moreover, it is $\frak{G}_{0}^{r+1}$-equivariant with
respect to the action of the group on $\mathcal{E}^{r}(V)$ defined above and
the natural left action of $\frak{G}_{0}^{r+1}$ on itself.
\end{lemma}

\textbf{Proof. }Let $(V,(x^{1},\ldots,x^{n}))$ be the standard chart of $V$.
Let us denote $f=\exp_{s}\circ u$, and $u=(X_{1},\ldots,X_{n})$, with
$X_{i}=\sum_{j}a_{i}^{j}\left(  \frac{\partial}{\partial x^{j}}\right)  _{0}$,
$i=1,\ldots,n$.

The map $t\mapsto f(tx)$ is the geodesic from $0\in V$ with initial speed
$u(x)$, i.e., the solution of the system of second order differential equations:%

\begin{equation}
\frac{d^{2}f^{k}(tx)}{dt^{2}}=-\sum_{j_{1},j_{2}}\Gamma_{j_{1}j_{2}}%
^{k}(f(tx))\frac{df^{j_{1}}(tx)}{dt}\frac{df^{j_{2}}(tx)}{dt},\hspace
{0.5cm}1\leq k\leq n, \label{geodesics}%
\end{equation}
with initial conditions $f(0)=0$, $\left.  \frac{df(tx)}{dt}\right|
_{0}=u(x)$. We can write this system in an equivalent way as follows:%

\begin{align*}
f(0)  &  =0\\
\sum_{i}x^{i}\frac{\partial f^{k}}{\partial x^{i}}(0)  &  =\sum_{i}x^{i}
a_{i}^{k},\hspace{0.5cm}1\leq k\leq n,\\
\sum_{i_{1},i_{2}}x^{i_{1}}x^{i_{2}}\frac{\partial^{2}f^{k}}{\partial
x^{i_{1}}\partial x^{i_{2}}}(tx)  &  =-\sum_{i_{1},i_{2},j_{1},j_{2}}x^{i_{1}
}x^{i_{2}}\Gamma_{j_{1}j_{2}}^{k}(f(tx))\frac{\partial f^{j_{1}}}{\partial
x^{i_{1}}}(tx)\frac{\partial f^{j_{2}}}{\partial x^{i_{2}}}(tx),~1\leq k\leq
n.
\end{align*}
From the second equation it follows that $\frac{\partial f^{k}}{\partial
x^{i}}(0)=a_{i}^{k}$, $1\leq i,k\leq n$. On the other hand, taking the
$(r-1)th$ order derivative ($r\geq1$) of the second equation and evaluating in
$t=0$ leads to the identity:
\[
\sum_{i_{1},\ldots,i_{r+1}}x^{i_{1}}\cdots x^{i_{r+1}}\frac{\partial
^{r+1}f^{k}}{\partial x^{i_{1}}\cdots\partial x^{i_{r+1}}}(0)=\sum
_{i_{1},\ldots i_{r+1}}x^{i_{1}}\cdots x^{i_{r+1}}R_{i_{1}\cdots i_{r+1}
},\hspace{0.5cm}1\leq k\leq n,
\]
where each $R_{i_{1}\cdots i_{r+1}}$ is a polynomial in the derivatives of the
Christoffel symbols, at $x=0,$ up to order $r-1,$ and the derivatives at $x=0$
of the components $f^{k}$ up to order $r$.

Last equation implies:%

\[
\frac{\partial^{r+1}f^{k}}{\partial x^{i_{1}}\cdots\partial x^{i_{r+1}}
}(0)=\frac{1}{(r+1)!}\sum_{\tau\in S_{r+1}}R_{i_{\tau(1)}\cdots i_{\tau(r+1)}
},\hspace{0.5cm}1\leq k,i_{1},\ldots,i_{r+1}\leq n.
\]
Since the map $j_{x}^{1}s\mapsto\nabla(s)(x)$ is smooth (see Remark
\ref{smoothness}), the derivatives up to order $r-1$ of the Christoffel
symbols $\Gamma_{ij}^{k}(x)$ are smooth functions of $j_{x}^{r}s$. Thus, the
$(r+1)th$ order derivatives at $x=0$ of $f$ are smooth functions of $j_{0}
^{r}s$ and $j_{0}^{r}f$.

Finally, using induction, it follows that $\mathrm{Exp}^{r}:\left(  j_{0}
^{r}s,u\right)  \mapsto j_{0}^{r+1}f$ is a well-defined smooth map.

In order to prove the $\frak{G}_{0}^{r+1}$-equivariance, let us point that,
due to the naturality of the assignment $s\mapsto\nabla(s)$, the equality
$f\circ\exp_{s}\circ u=\exp_{f\cdot s}\circ\bar{f}(u)$ holds for each
diffeomorphism $f:M\rightarrow M^{\prime}$ and each $G$-structure $s\in
\Gamma\left(  FM/G\right)  $. Therefore:%

\begin{align*}
\mathrm{Exp}^{r}\left(  \left(  j^{r+1}f\right)  \cdot\left(  j_{0}
^{r}s,u\right)  \right)   &  =\mathrm{Exp}^{r}\left(  j_{0}^{r}(f\cdot
s),\bar{f}(u)\right) \\
&  =j_{0}^{r+1}\left(  \exp_{f\cdot s}\circ\bar{f}(u)\right) \\
&  =j_{0}^{r+1}\left(  f\circ\exp_{s}\circ u\right) \\
&  =\left(  j^{r+1}f\right)  \cdot\mathrm{Exp}^{r}(j_{0}^{r}s,u),
\end{align*}
for each $j^{r+1}f\in\frak{G}_{0}^{r+1}$ and each $(j_{0}^{r}s,u)\in
\mathcal{E}^{r}(V)$.$\blacksquare$

Next, we define $\mathcal{E}^{r}(V)_{1}$ as the fiber over $1_{\frak{G}%
_{0}^{r+1}}\in\frak{G}_{0}^{r+1}$ of the map $\mathrm{Exp}^{r}:\mathcal{E}%
^{r}(V)\rightarrow\frak{G}_{0}^{r+1}$, i.e., $\mathcal{E}^{r}(V)_{1}:=\left(
\mathrm{Exp}^{r}\right)  ^{-1}\left(  1_{\frak{G}_{0}^{r+1}}\right)  $.

\begin{lemma}
$\mathcal{E}^{r}(V)_{1}$ is a smooth submanifold of $\mathcal{E}^{r}(V)$.
\end{lemma}

\textbf{Proof. }We prove first that, due to equivariance, the map
$\mathrm{Exp}^{r}$ is a surjective submersion, which implies that
$\mathcal{E}^{r}(V)_{1}$ is in fact a submanifold of $\mathcal{E}^{r}(V).$ It
is surjective because the image of each orbit in $\mathcal{E}^{r}(V)$ is an
orbit of left translations in the group, namely, $\frak{G}_{0}^{r+1}$. To see
that it is a submersion, note that if $(j_{0}^{r}s,u)\in\mathcal{E}^{r}(V)$
were a critical point of $\mathrm{Exp}^{r}$, then every element of its orbit
would also be a critical point, and, being surjective restricted to each
orbit, the map $\mathrm{Exp}^{r}$ would not have any regular value,
contradicting Sard's Theorem. Thus, each fiber of $\mathrm{Exp}^{r}$, in
particular $\mathcal{E}^{r}(V)_{1}$, is a smooth submanifold of $\mathcal{E}
^{r}(V)$.$\blacksquare$

On the other hand, the equivariance of $\mathrm{Exp}^{r}$ allows to define a
map $p^{r}:\mathcal{E}^{r}(V)\rightarrow\mathcal{E}^{r}(V)_{1}$ as
\[
p^{r}(j_{0}^{r}s,u)=\left(  \mathrm{Exp}^{r}(j_{0}^{r}s,u)\right)  ^{-1}
\cdot(j_{0}^{r}s,u).
\]
A straightforward calculation shows that this map is invariant under the right
action associated to the left action of $\frak{G}_{0}^{r+1}$ on $\mathcal{E}
^{r}(V),$ which is given by:
\[
\left(  j_{0}^{r}s,u\right)  \cdot\left(  j_{0}^{r+1}f\right)  =\left(
j_{0}^{r+1}f\right)  ^{-1}\cdot\left(  j_{0}^{r}s,u\right)  ,\hspace
{0.5cm}j_{0}^{r+1}f\in\frak{G}_{0}^{r+1},~\left(  j_{0}^{r}s,u\right)
\in\mathcal{E}^{r}(V).
\]
Indeed, the following lemma holds:

\begin{lemma}
$p^{r}:\mathcal{E}^{r}(V)\rightarrow\mathcal{E}^{r}(V)_{1}$ is a trivial
principal bundle with structure group $\frak{G}_{0}^{r+1}$.
\end{lemma}

\textbf{Proof.} The map $\Lambda^{r}:\mathcal{E}^{r}(V)\rightarrow
\mathcal{E}^{r}(V)_{1}\times\frak{G}_{0}^{r+1}$ defined by: $\Lambda
^{r}\left(  j_{0}^{r}s,u\right)  =\left(  p^{r}(j_{0}^{r}s,u),\left(
\mathrm{Exp}^{r}(j_{0}^{r}s,u)\right)  ^{-1}\right)  $ is a $\frak{G}
_{0}^{r+1}$-equivariant diffeomorphism. The inverse $\left(  \Lambda
^{r}\right)  ^{-1}:\mathcal{E}^{r}(V)_{1}\times\frak{G}_{0}^{r+1}
\rightarrow\mathcal{E}^{r}(V)$ is given by $\left(  \Lambda^{r}\right)
^{-1}\left(  \left(  j_{0}^{r}s,u\right)  ,j_{0}^{r+1}f\right)  =\left(
j_{0}^{r}s,u\right)  \cdot\left(  j_{0}^{r+1}f\right)  .$ This diffeomorphism
endows $\mathcal{E}^{r}(V)$ with the structure of a trivial principal
$\frak{G}_{0}^{r+1}$-bundle over $\mathcal{E}^{r}(V)_{1}$.$\blacksquare$

\begin{corollary}
The space $\mathbf{S}^{r}:=\mathcal{E}^{r}(V)/\frak{G}_{0}^{r+1}$ is a smooth
manifold canonically diffeomorphic to $\mathcal{E}^{r}(V)_{1}$.
\end{corollary}

\begin{remark}
\label{inversa}\emph{Explicitly, the diffeomorphism $\bar{p}^{r}%
:\mathbf{S}^{r}\rightarrow\mathcal{E}^{r}(V)_{1}$ is given by }%
\[
\overline{p}^{r}\left(  [(j_{0}^{r}s,u)]\right)  =p^{r}(j_{0}^{r}s,u).
\]
\emph{ Its inverse is the composed map: $\mathcal{E}^{r}(V)_{1}\hookrightarrow
\mathcal{E}^{r}(V)\rightarrow\mathcal{E}^{r}(V)/\frak{G}_{0}^{r+1}$ of the
inclusion (which is the map trivializing the bundle) with the natural projection.}
\end{remark}

\begin{remark}
\emph{Let us denote by $u^{0}:V\rightarrow T_{0}V$ the canonical frame. The
condition $j_{0}^{r+1}\left(  \exp_{s}\circ u\right)  =1_{\frak{G}_{0}^{r+1}}$
defining $\mathcal{E}^{r}(V)_{1}$ implies $D\left(  \exp_{s}\circ u\right)
(0)=\mathrm{id}_{V}$, where $D$ stands for the standard Fr\'{e}chet derivative
in $V$. Using that
\begin{align*}
D\left(  \exp_{s}\circ u\right)  (0) &  =D\left(  \exp_{s}\circ u^{0}\right)
(0)\circ D\left(  (u^{0})^{-1}\circ u\right)  (0)=\\
&  =D\left(  (u^{0})^{-1}\circ u\right)  (0)=(u^{0})^{-1}\circ u,
\end{align*}
we conclude that $u=u^{0}$, so that $\mathcal{E}^{r}(V)_{1}\subset J_{0}%
^{r}\left(  FV/G\right)  \times\left\{  u^{0}\right\}  $. Therefore from now
on, we will consider $\mathcal{E}^{r}(V)_{1}$ as a subspace of $J_{0}%
^{r}\left(  FV/G\right)  .$}
\end{remark}

\begin{remark}
\emph{Notice that an $r$-jet $j_{0}^{r}s\in J_{0}^{r}\left(  FV/G\right)  $
belongs to $\mathcal{E}^{r}(V)_{1}$ if and only if $s(0)=[u^{0}]$ and
\begin{equation}
\left(  \frac{d^{k-1}}{dt^{k-1}}\sum_{i_{1},i_{2}}x^{i_{1}}x^{i_{2}}%
\Gamma_{i_{1}i_{2}}^{\gamma}(tx)\right)  (0)=0,\hspace{1cm}1\leq k\leq
r,~1\leq\gamma\leq n,\label{r+1_NormalCoordinates}%
\end{equation}
for any $x$ in a neighborhood of $0\in V$, as follows from (\ref{geodesics}).}
\end{remark}

The action of $G$ on $\mathbf{S}^{r}$ can be translated now to $\mathcal{E}
^{r}(V)_{1}$.

\begin{lemma}
\label{B}The action of $G$ on $\mathbf{S}^{r}$ induces, by means of the
diffeomorphism $\mathbf{S}^{r}\cong\mathcal{E}^{r}(V)_{1}$, the following
action of $G$ on $\mathcal{E}^{r}(V)_{1}$:
\[
\left(  j_{0}^{r}s\right)  \cdot g=\left(  j_{0}^{r+1}g^{-1}\right)  \cdot
j_{0}^{r}s,\hspace{0.5cm}j_{0}^{r}s\in\mathcal{E}^{r}(V)_{1},~g\in
G\subset\frak{G}_{0}^{r+1}.
\]
\end{lemma}

\textbf{Proof.} For each $j_{0}^{r}s\in\mathcal{E}^{r}(V)_{1}$ and each $g\in
G$, we have, taking into account Remark \ref{inversa} we have that%

\begin{align*}
\left(  j_{0}^{r}s\right)  \cdot g  &  \doteqdot\overline{p}^{r}\left(
[(j_{0}^{r}s,u^{0})]\cdot g\right) \\
&  =\overline{p}^{r}\left(  [(j_{0}^{r}s,u^{0}\circ g)]\right) \\
&  =p^{r}\left(  j_{0}^{r}s,u^{0}\circ g\right) \\
&  =\left(  \mathrm{Exp}^{r}(j_{0}^{r}s,u^{0}\circ g)\right)  ^{-1}\cdot
(j_{0}^{r}s,u^{0}\circ g).
\end{align*}
From the fact that:%

\begin{align*}
\mathrm{Exp}^{r}(j_{0}^{r}s,u\circ g)  &  =j_{0}^{r+1}\left(  \exp_{s}\circ
u^{0}\circ g\right) \\
&  =j_{0}^{r+1}\left(  \exp_{s}\circ u^{0}\right)  \cdot j_{0}^{r+1}g\\
&  =j_{0}^{r+1}g,
\end{align*}
it follows the desired expression of $\left(  j_{0}^{r}s\right)  \cdot g:$%

\begin{align*}
\left(  j_{0}^{r}s\right)  \cdot g  &  =\left(  j_{0}^{r+1}g^{-1}\right)
\cdot(j_{0}^{r}s,u^{0}\circ g)\\
&  =\left(  \left(  j_{0}^{r+1}g^{-1}\right)  \cdot(j_{0}^{r}s),u^{0}\right)
\end{align*}
which is identified with $\left(  j_{0}^{r+1}g^{-1}\right)  \cdot(j_{0}
^{r}s).\blacksquare$

It is immediate that the projections $J_{0}^{r+k}\left(  FV/G\right)  \to
J_{0}^{r}\left(  FV/G\right)  $ induce projections $\mathbf{S}^{r+k}
\to\mathbf{S}^{r} $, $k\geq0$, which allow to define $\mathbf{S}^{\infty}=
\underset{\leftarrow}{\lim}\,\mathbf{S}^{r}$ as the projective limit of the
family $\{\mathbf{S}^{r}\}$. The action of $G$ induces an action on the limit
$\mathbf{S}^{\infty}$ in a natural way.

Thus, we have proved the following theorem of reduction of the group:

\begin{theorem}
The spaces $\mathbf{S}^{r}=\mathcal{E}^{r}(V)/\frak{G}_{0}^{r+1}$ have natural
smooth $G$-manifold structures, so that each moduli space $\frak{M}_{G}^{r}$
is canonically isomorphic to the quotient $\mathbf{S}^{r}/G$. Taking
projective limits yields a canonical isomorphism $\frak{M}_{G}^{\infty}%
\cong\mathbf{S}^{\infty}/G.$
\end{theorem}

\section{Embedding of the moduli of framed $r$-jets $\mathbf{S}^{r}$ in the
space of jets of moving frames $J_{0}^{r}(FV)$}

Let us denote by $J_{0}^{r}\left(  FV/G\right)  _{u^{0}}$ the manifold of
$r$-jets $j_{0}^{r}s\in J_{0}^{r}\left(  FV/G\right)  $ such that
$s(0)=[u^{0}]$. We will define an embedding of this manifold in $J_{0}
^{r}(FV)$. Then, from the diffeomorphism $\mathbf{S}^{r}\cong\mathcal{E}
^{r}(V)_{1}$ and the inclusion $\mathcal{E}^{r}(V)_{1}\hookrightarrow
J_{0}^{r}\left(  FV/G\right)  _{u^{0}}$, we will obtain an embedding of
$\mathbf{S}^{r}$ in $J_{0}^{r}(FV)$. This construction is a technical tool to
define the $G$-module structure on $\mathbf{S}^{r}$ in the next section.

We define $\Upsilon^{r}:J_{0}^{r}\left(  FV/G\right)  _{u^{0}}\rightarrow
J_{0}^{r}\left(  FV\right)  $ as follows. Given $j_{0}^{r}s\in J_{0}
^{r}\left(  FV/G\right)  _{u^{0}}$, let $s:\mathcal{U}\subset V\rightarrow
FV/G$ be a local representative of $j_{0}^{r}s$. We define a moving frame
$\sigma=\left(  X_{1},\ldots,X_{n}\right)  $ in a neighborhood of $0\in V$ as
the $\nabla(s)$-parallel transport of the frame $u^{0}$ along each ray
$x(t)=(t\cdot\lambda^{1},\ldots,t\cdot\lambda^{n})\in V$, and we set
$\Upsilon^{r}\left(  j_{0}^{r}s\right)  =j_{0}^{r}\sigma$.

\begin{lemma}
The assignment $\Upsilon^{r}:J_{0}^{r}\left(  FV/G\right)  _{u^{0}}\rightarrow
J_{0}^{r}\left(  FV\right)  $ is a well-defined embedding.
\end{lemma}

\textbf{Proof. }Let $s$ be a representative of $j_{0}^{r}s$ defined in a
neighborhood of $0\in V$, and denote by $\Gamma_{ij}^{k}$ the Christoffel
symbols of the connection $\nabla(s)$. Let $\sigma=\left(  X_{1},\ldots
,X_{n}\right)  $ be the moving frame defined above using the connection
$\nabla(s)$, and let $X_{i}=\sum_{k}\sigma_{ki}\left(  x\right)
\partial/\partial x^{k}$ be the expression of the vector field $X_{i}$ with
respect to the standard chart of $V$.

Since, by definition, each $X_{i}$ is parallel along the rays $x(t)=t\lambda$,
$\lambda=\left(  \lambda^{1},\ldots,\lambda^{n}\right)  $, the functions
$\sigma_{ki}$ satisfy the differential equations:%

\begin{align*}
\sum_{\beta}\left(  \dfrac{\partial\sigma_{\gamma i}}{\partial x^{\beta}%
}(t\lambda)+\sum_{\alpha}\Gamma_{\beta\alpha}^{\gamma}(t\lambda)\sigma_{\alpha
i}(t\lambda)\right)  \lambda^{\beta}  &  =0\\
\sigma_{ij}(0)  &  =\delta_{ij}%
\end{align*}
for any $\lambda\in V$, and any indices $1\leq\gamma,i,j\leq n$. Multiplying
these equations by $t,$ we obtain that $\sigma$ is characterized by the
initial condition $\sigma\left(  0\right)  =u^{0}$ and the equations:
\begin{equation}
\sum_{\beta}\left(  \frac{\partial\sigma_{\gamma i}}{\partial x^{\beta}}%
+\sum_{\alpha}\Gamma_{\beta\alpha}^{\gamma}\sigma_{\alpha i}\right)  x^{\beta
}=0,\text{ }1\leq i,\gamma\leq n. \label{Parallel}%
\end{equation}
Iterated derivation of (\ref{Parallel}) and evaluation in $x=0$ leads to:
\begin{equation}
\frac{\partial^{k}\sigma_{\gamma i}}{\partial x^{i_{1}}\cdots\partial
x^{i_{k}}}\left(  0\right)  =-\frac{1}{k}\underset{\left(  i_{1}\cdots
i_{k}\right)  }{\frak{S}}\left(  \frac{\partial^{k-1}}{\partial x^{i_{2}%
}\ldots\partial x^{i_{k}}}\sum_{\alpha}\Gamma_{i_{1}\alpha}^{\gamma}%
\sigma_{\alpha i}\right)  (0), \label{defUpsilon}%
\end{equation}
for any $k\geq1$ and any $1\leq\gamma,i,i_{1},\ldots,i_{r}\leq n$, were
$\frak{S}$ stands for the cyclic sum with respect to the corresponding indices.

For $r=1$ we have $\frac{\partial\sigma_{\gamma i}}{\partial x^{j}}
(0)=-\Gamma_{ji}^{\gamma}(0)$, so that $j_{0}^{1}\sigma=\Upsilon^{1}(j_{0}
^{1}s)$ depends smoothly on $j_{0}^{1}s.$ By induction, using
(\ref{defUpsilon}), it follows that $j_{0}^{r}\sigma=\Upsilon^{r}(j_{0}^{r}s)$
depends smoothly on $j_{0}^{r}s.$

Finally, it is obvious, by its very definition, that $\Upsilon^{r}$ takes its
values in the closed submanifold: $J_{0}^{r}\left(  FV\right)  _{u^{0}
}=\{j_{0}^{r}\sigma\in J_{0}^{r}\left(  FV\right)  \mid\sigma(0)=u^{0}\},$ and
that it is a section of the projection $J_{0}^{r}\left(  FV\right)  _{u^{0}
}\rightarrow J_{0}^{r}\left(  FV/G\right)  _{u^{0}}$ induced by
$q:FV\rightarrow FV/G$. Therefore, $\Upsilon^{r}$ is indeed an
embedding.$\blacksquare$

Next we define a right action of $G$ on $J_{0}^{r}\left(  FV\right)  _{u^{0}}$
by:
\[
\left(  j_{0}^{r}\sigma\right)  \cdot g=j_{0}^{r}\left(  R_{g}\circ\left(
g^{-1}\cdot\sigma\right)  \right)  ,\hspace{0.5cm}j_{0}^{r}\sigma\in J_{0}
^{r}\left(  FV\right)  _{u^{0}},\text{ }g\in G,
\]
where $R_{g}$ stands for right translations in $FV$. Besides, we have the
right action of $G$ on $J_{0}^{r}\left(  FV/G\right)  _{u^{0}}$:
\[
\left(  j_{0}^{r}s\right)  \cdot g=\left(  j_{0}^{r+1}g^{-1}\right)  \cdot
j_{0}^{r}s,\hspace{0.5cm}j_{0}^{r}s\in J_{0}^{r}\left(  FV/G\right)  _{u^{0}
},~g\in G.
\]

\begin{lemma}
The map $\Upsilon^{r}:J_{0}^{r}\left(  FV/G\right)  _{u^{0}}\rightarrow
J_{0}^{r}\left(  FV\right)  _{u^{0}}$ is $G$-equivariant with respect to the
right actions of $G$ previously defined.
\end{lemma}

\textbf{Proof. }Due to the naturality of the assignment $s\mapsto\nabla(s)$,
if $\sigma$ is the parallel moving frame associated to the $G$-structure $s$,
then $R_{g}\circ\left(  g^{-1}\cdot\sigma\right)  $ is the parallel moving
frame associated to $s\cdot g=\widetilde{g}^{-1}\circ s\circ g.$ From this
fact, the result follows easily.$\blacksquare$

Denoting $\mathcal{S}^{r}=\Upsilon^{r}\left(  \mathcal{E}^{r}(V)_{1}\right)
\subset J_{0}^{r}\left(  FV\right)  $ we have now the identifications between
$G$-manifolds $\mathbf{S}^{r}\cong\mathcal{E}^{r}(V)_{1}\cong\mathcal{S}^{r}$.

We derive now a characterization of the submanifolds $\mathcal{S}^{r}\subset
J_{0}^{r}\left(  FV\right)  $, which will be used in next section to define
inductively a $G$-module structure on each of them. First of all, we introduce
some definitions.

Let us denote $V^{r,p}_{q} := S^{r}(V^{\star})\otimes V^{\star\otimes
p}\otimes V^{\otimes q}$ with its natural $G$-module structure. Consider the
homomorphisms of $G$-modules $\mathrm{sym}^{r,p}_{q}:V^{r,p}_{q} \to
V^{r+1,p-1}_{q}$, $p \geq1$, given by symmetrization on the first $r+1$
covariant indices. We will denote $\bar{V}^{r,p}_{q} := \ker\left(
\mathrm{sym}^{r,p}_{q}\right)  $.

Let $\sigma=(X_{1},\ldots,X_{n})$, with $X_{i}=\sum_{j}\sigma_{ji}
\partial/\partial x^{j}$, $1\leq i\leq n$, be a representative of the $r$-jet
$j_{0}^{r}\sigma\in J_{0}^{r}\left(  FV\right)  $. Then, for each $g\in G$,
the derivatives of the components of the vectors in the representative
$R_{g}\circ\left(  g^{-1}\cdot\sigma\right)  $ of $\left(  j_{0}^{r}
\sigma\right)  \cdot g$ are given by $\ $
\[
\dfrac{\partial^{k}\bar{\sigma}_{rs}}{\partial x^{p_{1}}\ldots\partial
x^{p_{k}}}\left(  0\right)  =g^{r\alpha}\dfrac{\partial^{k}\sigma
_{\alpha\gamma}}{\partial x^{l_{1}}\ldots\partial x^{l_{k}}}\left(  0\right)
g_{l_{1}p_{1}}\cdots g_{l_{k}p_{k}}g_{\gamma s}.
\]
Thus, we can define an equivariant mapping $\Sigma^{(r)}:J_{0}^{r}\left(
FV\right)  \rightarrow V_{1}^{r,1}$ by:
\[
\Sigma^{(r)}\left(  j_{0}^{r}\sigma\right)  =\sum_{l_{1},\ldots,l_{r}}
\sum_{i,j}\frac{\partial^{r}\sigma_{ij}}{\partial x^{l_{1}}\ldots\partial
x^{l_{r}}}\left(  0\right)  v^{\star l_{1}}\odot\cdots\odot v^{\star l_{r}
}\otimes v^{\star j}\otimes v_{i},
\]
with $\{v_{1},\ldots,v_{n}\}$ being the canonical basis of $V$ and
$\{v^{\star1},\ldots,v^{\star n}\}$ its dual basis.

Next lemma is an immediate consequence of the above remarks:

\begin{lemma}
\label{DecompositionOfJetsOfFrames}The mapping $\left(  \pi_{r}^{r+1}%
,\Sigma^{(r+1)}\right)  :J_{0}^{r+1}\left(  FV\right)  \rightarrow J_{0}%
^{r}\left(  FV\right)  \times V_{1}^{r+1,1}$ is a $G$-equivariant diffeomorphism.
\end{lemma}

Now, let $\bar{\nabla}$ and $\nabla$ be the canonical connections adapted to
the $\{e\}$-structure $\sigma$ and the $G$-structure $s=q\circ\sigma$,
respectively. The connection $\bar{\nabla}$ is characterized by the equations
$\bar{\nabla}_{X_{i}}X_{j}=0$, $1\leq i,j\leq n$, so that the corresponding
Christoffel symbols are given by $\bar{\Gamma}_{\alpha\beta}^{\gamma}%
=-\sum_{i}\dfrac{\partial\sigma_{\gamma i}}{\partial x^{\alpha}}\sigma
^{i\beta}.$ We define the tensor field $F$ as the difference of both
connections, i.e., by $F(X,Y)=\nabla_{X}Y-\bar{\nabla}_{X}Y,\;X,Y\in
\frak{X}\left(  V\right)  $, and denote by $F_{\alpha\beta}^{\gamma}%
=\Gamma_{\alpha\beta}^{\gamma}-\bar{\Gamma}_{\alpha\beta}^{\gamma}$, its
components with respect to the canonical coordinates in $V$. From the
expression in local coordinates of the connection $\nabla$ given in Remark
\ref{forF}, it follows easily that these components have the form:
\begin{align}
F_{\alpha\beta}^{\gamma}  &  =\sum_{\lambda,\mu,\rho}\sigma^{\lambda\alpha
}\sigma^{\mu\beta}\sigma_{\gamma\rho}\left(  \sigma^{\star}\omega\right)
_{\mu}^{\rho}\left(  X_{\lambda}\right)  =\label{Fcomponents}\\
&  =\sum_{\lambda,\mu,\rho,i,j,k,h,\zeta}A_{\lambda\mu k}^{\rho ij}%
\sigma^{\lambda\alpha}\sigma^{\mu\beta}\sigma_{\gamma\rho}\sigma^{k\zeta
}\left(  \sigma_{hi}\frac{\partial\sigma_{\zeta j}}{\partial x^{h}}%
-\sigma_{hj}\frac{\partial\sigma_{\zeta i}}{\partial x^{h}}\right) \nonumber
\end{align}

Due to the naturality of $\bar{\nabla}$ and $\nabla$, for each $g\in G$, the
moving frame $R_{g}\circ\left(  g^{-1}\cdot\sigma\right)  $ has an associated
tensor field $\bar{F}$ whose components are related to those of $F$ by:
\[
\bar{F}_{ij}^{k}\left(  x\right)  =g^{k\gamma}F_{\alpha\beta}^{\gamma}\left(
g\cdot x\right)  g_{\alpha i}g_{\beta i}.
\]
Iterated derivation of these equations and evaluation at $x=0$ yields:
\[
\dfrac{\partial^{k}\bar{F}_{ij}^{h}}{\partial x^{p_{1}}\ldots\partial
x^{p_{k}}} (0) =g^{h\gamma}\dfrac{\partial^{k}F_{\alpha\beta}^{\gamma}}{
\partial x^{l_{1}}\ldots\partial x^{l_{k}}}\left(  0\right)  g_{l_{1}p_{1}
}\cdots g_{l_{k}p_{k}}g_{\alpha i}g_{\beta j}.
\]
So, we can also define an equivariant mapping $F^{\left(  r\right)  }
:J_{0}^{r+1}\left(  FV\right)  \rightarrow V^{r,2}_{1}$ as follows:
\[
F^{\left(  r\right)  }\left(  j_{0}^{r+1}\sigma\right)  =\sum_{l_{1}
,\ldots,l_{r}}\sum_{\alpha,\beta,\gamma}\frac{\partial^{r}F_{\alpha\beta
}^{\gamma}}{\partial x^{l_{1}}\ldots\partial x^{l_{r}}}\left(  0\right)
v^{\star l_{1}}\odot\cdots\odot v^{\star l_{r}}\otimes v^{\star\alpha} \otimes
v^{\star\beta} \otimes v_{\gamma}.
\]

\begin{lemma}
\label{decomposition} There exist a linear map $L^{(r)}:V_{1}^{r+1,1}%
\rightarrow V_{1}^{r,2}$ and a polynomial function $Q^{(r)}:J_{0}^{r}\left(
FV\right)  \rightarrow V_{1}^{r,2}$ in the standard coordinates of $J_{0}%
^{r}\left(  FV\right)  $ such that $F^{(r)}=L^{(r)}\circ\Sigma^{(r+1)}%
+Q^{(r)}\circ\pi_{r}^{r+1}$. Moreover, $L^{(r)}$ and $Q^{(r)}$ are $G$-equivariant.
\end{lemma}

\textbf{Proof. }Derivation of (\ref{Fcomponents}) and evaluation at $x=0$ show
that the derivatives of $F_{\alpha\beta}^{\gamma}$ at $x=0$ are of the form:
\[
\dfrac{\partial^{r}F_{\alpha\beta}^{\gamma}}{\partial x^{l_{1}}\ldots\partial
x^{l_{r}}}\left(  0\right)  =\sum_{i,j,k}A_{\alpha\beta k}^{\gamma ij}\left(
\frac{\partial^{r+1}\sigma_{kj}}{\partial x^{l_{1}}\ldots\partial x^{l_{r}
}\partial x^{i}}(0)-\frac{\partial^{r+1}\sigma_{ki}}{\partial x^{l_{1}}
\ldots\partial x^{l_{r}}\partial x^{j}}(0)\right)  +O_{r}
\]
where $O_{r}$ stands for a polynomial in the derivatives up to order $r$ of
$\sigma$ at $x=0$.

Accordingly, we define $L^{(r)}$ as the linear map sending each element of
$V_{1}^{r+1,1}$:
\[
\sum_{l_{1},\ldots,l_{r+1}}\sum_{j,k}t_{l_{1}\cdots l_{r+1}j}^{k}v^{\star
l_{1}}\odot\cdots\odot v^{\star l_{r+1}}\otimes v^{\star j}\otimes v_{k},
\]
to the element of $V_{1}^{r,2}$:
\[
\sum_{l_{1},\ldots,l_{r}}\sum_{\alpha,\beta,\gamma}\sum_{i,j,k}A_{\alpha\beta
k}^{\gamma ij}\left(  t_{l_{1}\cdots l_{r}ij}^{k}-t_{l_{1}\cdots l_{r}ji}%
^{k}\right)  v^{\star l_{1}}\odot\cdots\odot v^{\star l_{r}}\otimes
v^{\star\alpha}\otimes v^{\star\beta}\otimes v_{\gamma}.
\]
Thus, $L^{(r)}$ is the composition:
\[
V_{1}^{r+1,1}\overset{\delta^{r+1,1}}{\longrightarrow}S^{r}(V^{\star}%
)\otimes\bigwedge^{2}V^{\star}\otimes V\overset{1\otimes\left(  \delta
^{-1}\circ P_{\mathrm{im}\delta}\right)  }{\longrightarrow}S^{r}(V^{\star
})\otimes V^{\star}\otimes\frak{g}\subset V_{1}^{r,2},
\]
where the first arrow stands for the Spencer's operator, which is defined in
general as the map $\delta^{r+1,l-1}:S^{r+1}(V^{\star})\otimes\bigwedge
^{l-1}V^{\star}\otimes V\rightarrow S^{r}(V^{\star})\otimes\bigwedge
^{l}V^{\star}\otimes V$ given by:
\[
\left(  \delta^{r+1,l-1}t\right)  (u_{1},\ldots,u_{r},u_{1}^{\prime}%
,\ldots,u_{l}^{\prime})=\sum_{h=1}^{l}(-1)^{h+1}t(u_{1},\ldots,u_{r}%
,u_{h}^{\prime},u_{1}^{\prime},\ldots,\widehat{u^{\prime}}_{h},\ldots
,u_{l}^{\prime}),
\]
for each $t\in S^{r+1}(V^{\star})\otimes\bigwedge^{l-1}V^{\star}\otimes V$ and
any vectors $u_{1},\ldots,u_{r},u_{1}^{\prime},\ldots,u_{l}^{\prime}\in V$.

The difference $F^{(r)}-L^{(r)}\circ\Sigma^{(r+1)}:J_{0}^{r+1}\left(
FV\right)  \rightarrow V_{1}^{r,2}$ is a function depending polynomially in
the derivatives up to order $r$ of $\sigma$ at $0$, and hence there is a
well-defined polynomial function $Q^{(r)}:J_{0}^{r}\left(  FV\right)
\rightarrow V_{1}^{r,2}$ such that $F^{(r)}-L^{(r)}\circ\Sigma^{(r+1)}
=Q^{(r)}\circ\pi_{r}^{r+1}$.

Finally, it is clear that $L^{(r)}$ is $G$-equivariant, since it has been
defined as a composition of homomorphisms of $G$-modules. The equivariance of
the maps $F^{(r)}$, $L^{(r)}$, $\Sigma^{(r+1)}$ and $\pi_{r}^{r+1}$ implies
that $Q^{(r)}$ is also $G$-equivariant.$\blacksquare$

The desired characterization of the submanifold $\mathcal{S}^{r+1}$ is the following.

\begin{theorem}
\label{AnotherCharacterization}The submanifold $\mathcal{S}^{r+1}\subset
J_{0}^{r+1}\left(  FV\right)  $ can be characterized as the set of $\left(
r+1\right)  $-jets $j_{0}^{r+1}\sigma\in\left(  \pi_{r}^{r+1},\Sigma
^{(r+1)}\right)  ^{-1}\left(  \mathcal{S}^{r}\times\bar{V}_{1}^{r+1,1}\right)
$ such that
\[
\mathrm{sym}_{1}^{r,2}\left(  \left(  L^{(r)}\circ\Sigma^{(r+1)}\right)
\left(  j_{0}^{r+1}\sigma\right)  \right)  =-\mathrm{sym}_{1}^{r,2}\left(
\left(  Q^{(r)}\circ\pi_{r}^{r+1}\right)  \left(  j_{0}^{r+1}\sigma\right)
\right)  .
\]
\end{theorem}

\textbf{Proof.} By definition, an $r$-jet $j_{0}^{r}\sigma\in J_{0}^{r}\left(
FV\right)  $ belongs to $\mathcal{S}^{r}$ if and only if $j_{0}^{r}
\sigma=\Upsilon^{r}\left(  j_{0}^{r}s\right)  $ with $j_{0}^{r}s\in
\mathcal{E}^{r}(V)_{1}$. Since $\Upsilon^{r}$ is a section of the projection
$q^{r}:J_{0}^{r}\left(  FV\right)  _{u^{0}}\rightarrow J_{0}^{r}\left(
FV/G\right)  _{u^{0}}$, we have that $j_{0}^{r}s=q^{r}\left(  j_{0}^{r}
\sigma\right)  $. So that $j_{0}^{r}\sigma\in\mathcal{S}^{r}$ if and only if
$q^{r}\left(  j_{0}^{r}\sigma\right)  \in\mathcal{E}^{r}(V)_{1}$ and $\left(
\Upsilon^{r}\circ q^{r}\right)  \left(  j_{0}^{r}\sigma\right)  =j_{0}
^{r}\sigma$.

Let us first study the second condition, which means that, up to order $r,$
the moving frame $\sigma$ is obtained by the parallel displacement of the
canonical frame $u^{0}$ at $0\in V$ using the canonical connection of the
$G$-structure determined by $\sigma.$ For each $j_{0}^{r}\sigma\in J_{0}%
^{r}\left(  FV\right)  _{u^{0}}$, let $\sigma$ and $s$ be local
representatives of $j_{0}^{r}\sigma$ and $q^{r}\left(  j_{0}^{r}\sigma\right)
$, respectively, and let $\Gamma_{ij}^{k}$ be the Christoffel symbols of the
canonical connection $\nabla(s)$. Then, $\left(  \Upsilon^{r}\circ
q^{r}\right)  \left(  j_{0}^{r}\sigma\right)  =j_{0} ^{r}\sigma$ if and only
if $j_{0}^{r}\sigma$ satisfies equation (\ref{defUpsilon}) for any $1\leq
k\leq r$, $1\leq\gamma,i,i_{1},\ldots,i_{r}\leq n$.

Writing $\Gamma_{i_{1}\alpha}^{\gamma}=\bar{\Gamma}_{i_{1}\alpha}^{\gamma
}+F_{i_{1}\alpha}^{\gamma}=-\sum_{j}\sigma^{j\alpha}\dfrac{\partial
\sigma_{\gamma j}}{\partial x^{i_{1}}}+F_{i_{1}\alpha}^{\gamma}$, equation
(\ref{defUpsilon}) reads:
\[
\underset{\left(  i_{1}\cdots i_{k}\right)  }{\frak{S}}\left(  \frac
{\partial^{k-1}}{\partial x^{i_{2}}\ldots\partial x^{i_{k}}}\sum_{\alpha
}F_{i_{1}\alpha}^{\gamma}\sigma_{\alpha i}\right)  (0)=0.
\]
But this is equivalent to
\[
\sum_{i_{1},\ldots,i_{k}}x^{i_{1}}\cdots x^{i_{k}}\left(  \frac{\partial
^{k-1}}{\partial x^{i_{2}}\ldots\partial x^{i_{k}}}\sum_{\alpha}F_{i_{1}
\alpha}^{\gamma}\sigma_{\alpha i}\right)  (0)=0,
\]
and hence to the equations
\[
\left(  \frac{d^{k-1}}{dt^{k-1}}\sum_{\alpha}\sigma_{\alpha i}(tx)\sum_{\beta
}x^{\beta}F_{\beta\alpha}^{\gamma}(tx)\right)  (0)=0,
\]
for any $1\leq k\leq r$, $1\leq\gamma,i,\leq n$ and any $x$ in a neighborhood
of $0\in V$. Using induction, these equations are seen to be equivalent to the following:%

\[
\left(  \frac{d^{k-1}}{dt^{k-1}}\sum_{\beta}x^{\beta}F_{\beta i}^{\gamma
}(tx)\right)  (0)=0,
\]
for any $1\leq k\leq r$, $1\leq\gamma,i,\leq n$ and any $x$, which can be
written in terms of the partial derivatives at $0\in V$ of $F$ (and hence of
$\sigma$) as:%

\begin{equation}
\underset{\left(  i_{1}\cdots i_{k}\right)  }{\frak{S}}\dfrac{\partial
^{k-1}F_{i_{1}\alpha}^{\gamma}}{\partial x^{i_{2}}\cdots\partial x^{i_{k}}%
}(0)=0, \label{SecondEquations}%
\end{equation}
for each $1\leq k\leq r$, and for any indices $1\leq\alpha,\gamma,i_{1}%
,\ldots,i_{k}\leq n$.

Let us now find the equations provided by the condition $q^{r}\left(
j_{0}^{r}\sigma\right)  \in\mathcal{E}^{r}(V)_{1},$ which means that, up to
order $r+1,$ the normal coordinates of the connection attached to the
$G$-structure $G\cdot\sigma$ are just the canonical coordinates of $V.$
Setting $\Gamma_{i_{1}i_{2}}^{\gamma}=\bar{\Gamma}_{i_{1}i_{2}}^{\gamma
}+F_{i_{1}i_{2}}^{\gamma}=-\sum_{j}\sigma^{ji_{2}}\dfrac{\partial
\sigma_{\gamma j}}{\partial x_{1}^{i}}+F_{i_{1}i_{2}}^{\gamma}=\sum_{j}%
\sigma_{\gamma j}\dfrac{\partial\sigma^{ji_{2}}}{\partial x_{1}^{i}}%
+F_{i_{1}i_{2}}^{\gamma}$ in (\ref{r+1_NormalCoordinates}) and using
(\ref{SecondEquations}), we obtain that $q^{r}\left(  j_{0}^{r}\sigma\right)
\in\mathcal{E}^{r}(V)_{1}$ if and only if:
\[
\left(  \frac{d^{k-1}}{dt^{k-1}}\sum_{j}\sigma_{\gamma j}(tx)\frac{d}%
{dt}\left(  \sum_{i_{2}}x^{i_{2}}\sigma^{ji_{2}}(tx)\right)  \right)  (0)=0,
\]
for any $1\leq k\leq r$, $1\leq\gamma\leq n$, and any $x$ in a neighborhood of
$0\in V$. By induction, it follows again that these equations are equivalent
to:
\[
\left(  \frac{d^{k}}{dt^{k}}\sum_{\beta}x^{\beta}\sigma^{\gamma\beta
}(tx)\right)  (0)=0,
\]
for any $1\leq k\leq r$, $1\leq\gamma\leq n$, and any $x$, which can be
written as:
\[
\underset{\left(  i_{1}\cdots i_{k+1}\right)  }{\frak{S}}\dfrac{\partial
^{k}\sigma^{\gamma i_{1}}}{\partial x^{i_{2}}\cdots\partial x^{i_{k+1}}}(0)=0
\]
or, equivalently, as:%

\begin{equation}
\underset{\left(  i_{1}\cdots i_{k+1}\right)  }{\frak{S}}\dfrac{\partial
^{k}\sigma_{\gamma i_{1}}}{\partial x^{i_{2}}\cdots\partial x^{i_{k+1}}}(0)=0
\label{FirstEquations}%
\end{equation}
for each $1\leq k\leq r$, and for any indices $1\leq\gamma,i_{1}%
,\ldots,i_{k+1}\leq n$. Finally, an $r$-jet $j_{0}^{r}\sigma\in J_{0}%
^{r}\left(  FV\right)  $ belongs to $\mathcal{S}^{r}$ if and only if it
satisfies $\sigma_{\alpha\beta}(0)=\delta_{\alpha\beta}$, together with
equations (\ref{FirstEquations}) and (\ref{SecondEquations}).

For the inductive definition of the submanifolds $\mathcal{S}^{r}$ appearing
in the statement of the theorem, notice that an $\left(  r+1\right)  $-jet
$j_{0}^{r+1}\sigma\in J_{0}^{r+1}\left(  FV\right)  $ lies in $\mathcal{S}%
^{r+1}$ if and only if $\pi_{r}^{r+1}\left(  j_{0}^{r+1}\sigma\right)
\in\mathcal{S}^{r}$ and equations (\ref{FirstEquations}) and
(\ref{SecondEquations}) are satisfied at the top level, i.e.: $\ $
\begin{equation}
\underset{\left(  i_{1}\cdots i_{r+2}\right)  }{\frak{S}}\dfrac{\partial
^{r+1}\sigma_{\alpha i_{1}}}{\partial x^{i_{2}}\cdots\partial x^{i_{r+2}}%
}(0)=0, \label{FirstTopEquations}%
\end{equation}%
\begin{equation}
\underset{\left(  i_{1}\cdots i_{r+1}\right)  }{\frak{S}}\dfrac{\partial
^{r}F_{i_{1}\alpha}^{\gamma}}{\partial x^{i_{2}}\cdots\partial x^{i_{r+1}}%
}\left(  0\right)  =0, \label{SecondTopEquations}%
\end{equation}
for any indices $1\leq\alpha,\gamma,i_{1},\ldots,i_{r+2}\leq n$.

Using the diffeomorphism given in Lemma \ref{DecompositionOfJetsOfFrames} and
the definition of the operators $F^{(r)}$, and taking into account that
\begin{align*}
&  \mathrm{sym}_{q}^{r,p}\left(  v^{\star i_{1}}\odot\cdots\odot v^{\star
i_{r}}\otimes v^{\star i_{r+1}}\otimes\cdots\otimes v^{\star i_{r+p}}\otimes
v_{j_{1}}\otimes\cdots\otimes v_{j_{q}}\right) \\
&  =v^{\star i_{1}}\odot\cdots\odot v^{\star i_{r}}\odot v^{\star i_{r+1}%
}\otimes\cdots\otimes v^{\star i_{r+p}}\otimes v_{j_{1}}\otimes\cdots\otimes
v_{j_{q}}\\
&  =\frac{1}{r+1}\underset{\left(  i_{1}\cdots i_{r+1}\right)  }{\frak{S}%
}v^{\star i_{1}}\odot\cdots\odot v^{\star i_{r}}\otimes v^{\star i_{r+1}%
}\otimes\cdots\otimes v^{\star i_{r+p}}\otimes v_{j_{1}}\otimes\cdots\otimes
v_{j_{q}},
\end{align*}
we see that equations (\ref{FirstTopEquations}) and (\ref{SecondTopEquations})
are equivalent to $\mathrm{sym}_{1}^{r+1,1}\left(  \Sigma^{(r+1)}\left(
j_{0}^{r+1}\sigma\right)  \right)  =0$ and $\mathrm{sym}_{1}^{r,2}\left(
F^{(r)}\left(  j_{0}^{r+1}\sigma\right)  \right)  $ $=0,$ so each submanifold
$\mathcal{S}^{r+1}\subset J_{0}^{r+1}\left(  FV\right)  $ is characterized by
these conditions together with $\pi_{r}^{r+1}\left(  j_{0}^{r+1}\sigma\right)
\in\mathcal{S}^{r}$, . That is:
\[
\mathcal{S}^{r+1}=\left\{  j_{0}^{r+1}\sigma\in\left(  \pi_{r}^{r+1}%
,\Sigma^{(r+1)}\right)  ^{-1}\left(  \mathcal{S}^{r}\times\bar{V}_{1}%
^{r+1,1}\right)  \mid\mathrm{sym}_{1}^{r,2}\left(  F^{(r)}\left(  j_{0}%
^{r+1}\sigma\right)  \right)  =0\right\}
\]

Finally, using the decomposition of $F^{(r)}=L^{(r)}\circ\Sigma^{(r+1)}
+Q^{(r)}\circ\pi_{r}^{r+1}$ given in Lemma \ref{decomposition}, we obtain the
desired result.$\blacksquare$

\section{\strut Definition of the $G$-module structure}

In this section, we will define a $G$-module structure on each $\mathcal{S}
^{r}$, inductively on $r$, such that the projections are homomorphisms of
$G$-modules. The corresponding $G$-manifolds $\mathbf{S}^{r}=\mathcal{E}
^{r}(V)/\frak{G}_{0}^{r+1}$ will inherit $G$-module structures, and the
natural projections between them will also be homomorphisms of $G$-modules.
Therefore we will obtain a $G$-module structure on the projective limit
$\mathbf{S}^{\infty}=\underset{\leftarrow}{\lim\,}\mathbf{S}^{r}$.

A main tool to define a $G$-module structure on each $\mathcal{S}^{r}$ (and
therefore on each $\mathbf{S}^{r}$) is next lemma.

\begin{lemma}
\label{section} The projections $\pi_{r}^{r+1}:\mathcal{S}^{r+1}%
\rightarrow\mathcal{S}^{r}$ admit smooth $G$-equivariant sections $\mu
_{r+1}^{r}:\mathcal{S}^{r}\rightarrow\mathcal{S}^{r+1}$.
\end{lemma}

\textbf{Proof. }Let us denote by $\bar{L}^{(r)}:\bar{V}_{1}^{r+1,1}\rightarrow
V_{1}^{r,2}$ the restriction of $L^{(r)}$ to $\bar{V}_{1}^{r+1,1}$, and
$W^{r+1}:=\ker\left(  \mathrm{sym}_{1}^{r,2}\circ\bar{L}^{(r)}\right)  $. Let
us assume for the moment that we can choose a $G$-submodule $Z^{r+1}%
\subset\bar{V}_{1}^{r+1,1}$ such that $\bar{V}_{1}^{r+1,1}=W^{r+1}\oplus
Z^{r+1}$. Then, the restriction $\left.  \left(  \mathrm{sym}_{1}^{r,2}%
\circ\bar{L}^{(r)}\right)  \right|  _{Z^{r+1}}:Z^{r+1}\rightarrow
\mathrm{im}\left(  \mathrm{sym}_{1}^{r,2}\circ\bar{L}^{(r)}\right)  $ is an
isomorphism of $G$-modules. Using that the projections $\pi_{r}^{r+1}%
:\mathcal{S}^{r+1}\rightarrow\mathcal{S}^{r}$ are onto, together with the fact
that, by Theorem \ref{AnotherCharacterization} we have that $\mathrm{sym}%
_{1}^{r,2}\circ L^{(r)}\circ\Sigma^{(r+1)}=-\mathrm{sym}_{1}^{r,2}\circ
Q^{(r)}\circ\pi_{r}^{r+1}$ on $\mathcal{S}^{r+1}$ and $\Sigma^{(r+1)}\left(
\mathcal{S}^{r+1}\right)  \subset\bar{V}_{1}^{r+1,1}$, we obtain that
\begin{align*}
\left(  -\mathrm{sym}_{1}^{r,2}\circ Q^{(r)}\right)  \left(  \mathcal{S}%
^{r}\right)   &  =\left(  -\mathrm{sym}_{1}^{r,2}\circ Q^{(r)}\circ\pi
_{r}^{r+1}\right)  \left(  \mathcal{S}^{r+1}\right)  \\
&  =\left(  \mathrm{sym}_{1}^{r,2}\circ L^{(r)}\circ\Sigma^{(r+1)}\right)
\left(  \mathcal{S}^{r+1}\right)  \subset\mathrm{im}\left(  \mathrm{sym}%
_{1}^{r,2}\circ\bar{L}^{(r)}\right)  .
\end{align*}
Thus, for each $j_{0}^{r}\sigma\in\mathcal{S}^{r}$ there exists a unique
$v^{r+1}\in Z^{r+1}$ such that
\[
\left(  \mathrm{sym}_{1}^{r,2}\circ\bar{L}^{(r)}\right)  \left(
v^{r+1}\right)  =-\mathrm{sym}_{1}^{r,2}\left(  Q^{\left(  r\right)  }\left(
j_{0}^{r}\sigma\right)  \right)  .
\]
Now, using the isomorphism of Lemma \ref{DecompositionOfJetsOfFrames}, we set
\[
\mu_{r+1}^{r}\left(  j_{0}^{r}\sigma\right)  =\left(  \pi_{r}^{r+1}%
\times\Sigma^{(r+1)}\right)  ^{-1}\left(  j_{0}^{r}\sigma,v^{r+1}\right)  .
\]

It is clear that this defines a section $\mu_{r+1}^{r}:\mathcal{S}%
^{r}\rightarrow\mathcal{S}^{r+1}$ of the projection $\pi_{r}^{r+1}$. Moreover,
from the $G$-equivariance of $\pi_{r}^{r+1}$, $\Sigma^{(r+1)}$, \textrm{sym}
$_{1}^{r,2}$, $Q^{(r)}$, $L^{(r)}$ and the fact that $Z^{r+1}$ is a
$G$-submodule, it follows that $\mu_{r+1}^{r}$is $G$-equivariant too.

It only remains to prove that it is indeed possible to choose a $G$-submodule
$Z^{r+1}\subset\bar{V}_{1}^{r+1,1}$ such that $\bar{V}_{1}^{r+1,1}%
=W^{r+1}\oplus Z^{r+1}.$ Note first that
\begin{equation}
V_{1}^{r+1,1}=\ker\left(  \mathrm{sym}_{1}^{r+1,1}\right)  \oplus\left(
S^{r+2}\left(  V^{\star}\right)  \otimes V\right)  =\bar{V}_{1}^{r+1,1}%
\oplus\ker\delta^{r+1,1}.\label{descomposicion}%
\end{equation}
Therefore, $\left.  \delta\right|  _{\bar{V}_{1}^{r+1,1}}:\bar{V}_{1}%
^{r+1,1}\longrightarrow\mathrm{im}\delta^{r+1,1}$ is an isomorphism of
$G$-modules. It is easy to check that its inverse is given by $\left(
\frac{r+1}{r+2}\right)  \mathrm{sym}_{1}^{r,2}:\mathrm{im}\delta
^{r+1,1}\longrightarrow\bar{V}_{1}^{r+1,1}.$ 

Hence, it will enough to find a supplementary $G$-submodule $\hat{Z}^{r+1}$
of
\[
\ker\left(  \left.  \mathrm{sym}_{1}^{r,2}\circ\left(  1\otimes\delta
^{-1}\circ P_{\mathrm{im}\delta}\right)  \right|  _{\mathrm{im}\delta^{r+1,1}%
}\right)
\]
in $\mathrm{im}\delta^{r+1,1}.$ In order to do this, let us first prove that,
if we consider $S^{r+1}\left(  V^{\star}\right)  \otimes\frak{g}$ as a
$G$-submodule of $V_{1}^{r+1,1},$ then
\[
\left.  \mathrm{sym}_{1}^{r,2}\circ L^{\left(  r\right)  }\right|
_{S^{r+1}\left(  V^{\star}\right)  \otimes\frak{g}}=1_{S^{r+1}\left(
V^{\star}\right)  \otimes\frak{g}}.
\]
So let $\tau\in S^{r+1}\left(  V^{\star}\right)  \otimes\frak{g.}$ Given
vectors $u_{1},\ldots,u_{r+1},u_{r+2}\in V,$ we can consider $\tau\left(
u_{1},\ldots,u_{r+1}\right)  \in\frak{g,}$ and so $\tau\left(  u_{1}%
,\ldots,u_{r+1}\right)  u_{r+2}\in V.$ We have that%
\begin{align*}
\left(  \delta^{r+1,1}\tau\right)  \left(  u_{1},\ldots,u_{r+1},u_{r+2}%
\right)   &  =\tau\left(  u_{1},\ldots,u_{r+1}\right)  u_{r+2}-\tau\left(
u_{1},\ldots,u_{r+2}\right)  u_{r+1}\\
&  =\delta\left(  \tau\left(  u_{1},\ldots,u_{r},\cdot\right)  \right)
\left(  u_{r+1},u_{r+2}\right)  \\
&  =\left[  \left(  1_{S^{r}\left(  V^{\star}\right)  }\otimes\delta\right)
\tau\right]  \left(  u_{1},\ldots,u_{r},u_{r+1},u_{r+2}\right)  .
\end{align*}
Therefore, $\delta^{r+1,1}\tau=\left(  1_{S^{r}\left(  V^{\star}\right)
}\otimes\delta\right)  \tau,$ whence%
\begin{align}
\left(  \mathrm{sym}_{1}^{r,2}\circ L^{\left(  r\right)  }\right)  \tau &
=\left(  \mathrm{sym}_{1}^{r,2}\circ\left(  1_{S^{r}\left(  V^{\star}\right)
}\otimes\delta^{-1}\circ P_{\mathrm{im}\delta}\right)  \circ\delta
^{r+1,1}\right)  \tau\label{identity}\\
&  =\mathrm{sym}_{1}^{r,2}\tau=\tau.\nonumber
\end{align}
Thus,
\[
\mathrm{im}\left(  \mathrm{sym}_{1}^{r,2}\circ L^{\left(  r\right)  }\right)
=S^{r+1}\left(  V^{\star}\right)  \otimes\frak{g}%
\]
and, according to (\ref{descomposicion}), we also have that
\[
\mathrm{im}\left(  \mathrm{sym}_{1}^{r,2}\circ\bar{L}^{\left(  r\right)
}\right)  =\mathrm{im}\left(  \mathrm{sym}_{1}^{r,2}\circ L^{\left(  r\right)
}\right)  =S^{r+1}\left(  V^{\star}\right)  \otimes\frak{g.}%
\]
Taking into account that, by (\ref{identity}), $\left.  \delta^{r+1,1}\right|
_{S^{r+1}\left(  V^{\star}\right)  \otimes\frak{g}}$ is injective, we conclude
that%
\begin{align*}
&  \dim\ker\left(  \left.  \mathrm{sym}_{1}^{r,2}\circ\left(  1\otimes
\delta^{-1}\circ P_{\mathrm{im}\delta}\right)  \right|  _{\mathrm{im}%
\delta^{r+1,1}}\right)  \\
&  =\dim\left(  \mathrm{im}\delta^{r+1,1}\right)  -\dim\left(  S^{r+1}\left(
V^{\star}\right)  \otimes\frak{g}\right)  \\
&  =\dim\left(  \mathrm{im}\delta^{r+1,1}\right)  -\dim\left(  \delta
^{r+1,1}\left(  S^{r+1}\left(  V^{\star}\right)  \otimes\frak{g}\right)
\right)  .
\end{align*}
Let us define $\hat{Z}^{r+1}=\delta^{r+1,1}\left(  S^{r+1}\left(  V^{\star
}\right)  \otimes\frak{g}\right)  .$ Obviously it is a $G$-submodule of
$\mathrm{im}\delta^{r+1,1}$ of the required dimension. Moreover,
(\ref{identity}) implies that%
\[
\hat{Z}^{r+1}\cap\ker\left(  \left.  \mathrm{sym}_{1}^{r,2}\circ\left(
1\otimes\delta^{-1}\circ P_{\mathrm{im}\delta}\right)  \right|  _{\mathrm{im}%
\delta^{r+1,1}}\right)  =\left\{  0\right\}  .
\]
Hence, $\hat{Z}^{r+1}$ is the desired supplementary.$\blacksquare$

\noindent\textbf{Proof of Theorem \ref{Theorem}:} The proof will be by
induction on $r.$ The case $r=0$ is trivial: since $\mathcal{S}^{0}$ has a
unique element, determined by the conditions $\sigma_{ij}\left(  0\right)
=\delta_{ij}$, it admits the trivial $G$-module structure.

Now, assume that a $G$-module structure has been already defined on
$\mathcal{S}^{r}$ and let us define one on $\mathcal{S}^{r+1}$.

Let us consider the $G$-equivariant map $f^{r}=\Sigma^{(r+1)}\circ\mu
_{r+1}^{r}:\mathcal{S}^{r} \to\bar{V}^{r+1,1}_{1}$, and define a
$G$-equivariant bijection $\psi^{r}:\mathcal{S}^{r}\times\bar{V}^{r+1,1}
_{1}\to\mathcal{S}^{r}\times\bar{V}^{r+1,1}_{1}$ by $\psi^{r}\left(  j_{0}
^{r}\sigma, v^{r+1}\right)  = \left(  j_{0}^{r}\sigma, v^{r+1} - f^{r}\left(
j_{0}^{r}\sigma\right)  \right)  $, where the action of $G$ on both sides is
just the diagonal one.

After Lemma \ref{section}, we see that the image $\left(  \pi_{r}^{r+1}
\times\Sigma^{(r+1)}\right)  \left(  \mathcal{S}^{r+1}\right)  $ can be
characterized as the set of pairs $\left(  j_{0}^{r}\sigma,v^{r+1}\right)
\in\mathcal{S}^{r}\times\bar{V}_{1}^{r+1,1}$ satisfying:
\[
\left(  \mathrm{sym}_{1}^{r,2}\circ\bar{L}^{(r)}\right)  (v^{r+1})=\left(
\mathrm{sym}_{1}^{r,2}\circ\bar{L}^{(r)}\right)  \left(  f^{r}\left(
j_{0}^{r}\sigma\right)  \right)
\]
(both members being equal to $-\mathrm{sym}_{1}^{r,2}\left(  Q^{\left(
r\right)  }\left(  j_{0}^{r}\sigma\right)  \right)  $). But this set is just
the preimage by $\psi^{r}$ of the $G$-submodule $\mathcal{S}^{r} \oplus
W^{r+1} \subset\mathcal{S}^{r}\oplus\bar{V}_{1}^{r+1,1}$.

Now, we define the $G$-module structure on $\mathcal{S}^{r+1}$ as the one that
makes linear the bijection:
\begin{equation}
\phi^{r}:=\psi^{r}\circ\left(  \pi_{r}^{r+1}\times\Sigma^{(r+1)}\right)
:\mathcal{S}^{r+1}\rightarrow\mathcal{S}^{r}\oplus W^{r+1},
\end{equation}
where we are considering the direct sum $G$-module structure on the right hand
side .

The following diagram
\[%
\begin{array}
[c]{rcl}%
\mathcal{S}^{r+1} & \overset{\phi^{r}}{\longrightarrow} & \mathcal{S}%
^{r}\oplus W^{r+1}\\
{\scriptsize \pi}_{r}^{r+1}\searrow &  & \swarrow{\scriptsize pr}_{1}\\
& \mathcal{S}^{r} &
\end{array}
\]
where $\mathrm{pr}_{1}$ denotes the projection onto the first summand, is
commutative, so that the projections $\pi_{r}^{r+1}$ are linear. Using
induction, as well as the fact that $\phi^{r}$ is $G$-equivariant, we conclude
that the action of $G$ on each $\mathcal{S}^{r}$ is also linear. Moreover, the
smoothness of the maps $f^{r}$ implies the compatibility of the $G$-module
structure with the manifold structure on each $\mathcal{S}^{r}$.$\blacksquare$

Finally, the resulting $G$-module structure on each $G$-manifold
$\mathbf{S}^{r}$ does not depend either on the supplementary subspace $W$ used
to define the canonical connections or on the sections $\mu_{r+1}^{r}$ between
the corresponding spaces $\mathcal{S}^{r}$. In fact, changing any of them
would define isomorphic $G$-module structures, as follows from next general lemma:

\begin{lemma}
\label{uniqueness} Let $M$ and $N$ be two smooth $n$-dimensional manifolds
endowed with smooth $G$-module structures. If $M$ and $N$ are diffeomorphic as
$G$-manifolds, then they also are isomorphic as $G$-modules.
\end{lemma}

\textbf{Proof.} Let $f:M\rightarrow N$ be a $G$-equivariant diffeomorphism. By
composing, if necessary, with a translation in $N$ we can assume that
$f(0)=0$. Due to the smoothness of the module structures, any linear
isomorphisms $\varphi_{1}:M\rightarrow\mathbb{R}^{n}$, $\varphi_{2}
:N\rightarrow\mathbb{R}^{n}$ (onto $\mathbb{R}^{n}$ with its standard linear
structure) define global charts of $M$ and $N$. The map $\varphi_{2}\circ
f\circ\varphi_{1}^{-1}:\mathbb{R}^{n}\rightarrow\mathbb{R}^{n}$ is a
diffeomorphism leaving $0$ fixed, and it is equivariant with respect to the
linear actions of $G$ on $\mathbb{R}^{n}$ induced by $\varphi_{1}$ and
$\varphi_{2}$. Then, the Fr\'{e}chet derivative $D\left(  \varphi_{2}\circ
f\circ\varphi_{1}^{-1}\right)  (0):\mathbb{R}^{n}\rightarrow\mathbb{R}^{n}$ is
a linear $G$-equivariant isomorphism, and so is the composed map: $\varphi
_{2}^{-1}\circ D\left(  \varphi_{2}\circ f\circ\varphi_{1}^{-1}\right)
(0)\circ\varphi_{1}:M\rightarrow N.$$\blacksquare$

By definition, the $G$-module structure on $\mathcal{S}^{r+1}$ is isomorphic
to $\mathcal{S}^{r}\oplus W^{r+1}$. This fact (together with the obvious
identity $\mathcal{S}^{0}=\{0\}$) yields an isomorphism of $G$-modules:
\[
\mathbf{S}^{r}\cong W^{1}\oplus W^{2}\oplus\cdots\oplus W^{r},
\]
and the $G$-module structure of $\mathbf{S}^{r}$ is determined by that of the
spaces $W^{k}$, $k=1,\ldots, r$.

An element $t\in W^{r+1}$ lies in $\bar{V}_{1}^{r+1,1}:=\ker\left(
\mathrm{sym}_{1}^{r+1,1}\right)  ,$ whence its components must satisfy the
equations
\begin{equation}
\underset{\left(  i_{1}\cdots i_{r+1}i_{r+2}\right)  }{\frak{S}}t_{i_{1}\cdots
i_{r+1}i_{r+2}}^{k}~=0,\hspace{1cm}1\leq k,i_{1},\ldots,i_{r+2}\leq n.
\label{W1}%
\end{equation}
Denoting by $(L^{(r)}t)_{i_{1}\cdots i_{r+1}i_{r+2}}^{k}$ the components of
$L^{(r)}t$, it follows that an element $t\in\bar{V}_{1}^{r+1,1}$ lies in
$W^{r+1}$ if and only if the following equations are also satisfied:
\begin{equation}
\underset{\left(  i_{1}\cdots i_{r+1}\right)  }{\frak{S}}(L^{(r)}%
t)_{i_{1}\cdots i_{r+1}i_{r+2}}^{k}=0,\hspace{1cm}1\leq k,i_{1},\ldots
,i_{r+2}\leq n. \label{W2}%
\end{equation}
Therefore, $W^{r+1}$ is the $G$-submodule of $V_{1}^{r+1,1}=S^{r+1}(V^{\star
})\otimes V^{\star}\otimes V$ characterized by equations (\ref{W1}) and
(\ref{W2}).

Notice that, in particular, $W^{1}=\ker\left(  \mathrm{sym}_{1}^{0,2}\circ
\bar{L}^{(0)}\right)  .$ By definition, $\mathrm{sym}_{1}^{0,2}$ is the
identity in $V^{\star}\otimes V^{\star}\otimes V.$ The mapping $\bar{L}%
^{(0)}:\bar{V}_{1}^{1,1}\rightarrow V_{1}^{0,2}$ is defined on $\bar{V}%
_{1}^{1,1}=\bigwedge^{2}V^{\star}\otimes V$ as
\[
\bigwedge\nolimits^{2}V^{\star}\otimes V\overset{P_{\mathrm{im}\delta}%
}{\longrightarrow}\delta\left(  V^{\star}\otimes\frak{g}\right)
\overset{\delta^{-1}}{\longrightarrow}V^{\star}\otimes\frak{g.}%
\]
Thus $W^{1}=\ker P_{\mathrm{im}\delta}=W.$

%
%
%
%

Next, we will describe the $G$-modules $W^{r}$ for some different choices of
$G$.

\begin{ejemplo}
[$\{e\}$-structures]\emph{For complete parallelisms, $W^{r+1}$ is
characterized by equation (\ref{W1}). Equation (\ref{W2}) is empty because the
tensor $F$, and hence the operator $L^{(r)}$, are zero. Therefore,
$W^{r+1}=\bar{V}_{1}^{r+1,1}$. }
\end{ejemplo}

\begin{ejemplo}
[$O(n)$-structures]\label{detailed} \emph{We have seen (Example
\ref{Levi-Civita}) that, for $G=O(n)$, the map $\delta$ is an isomorphism,
with inverse $\delta^{-1}:\bigwedge^{2}V^{\star}\otimes V\rightarrow V^{\star
}\otimes\frak{g}$ given, in components, by:
\[
\left(  \delta^{-1}(T)\right)  _{ij}^{k}=\frac{1}{2}\left(  T_{ij}^{k}%
+T_{ki}^{j}+T_{kj}^{i}\right)
\]
}

\emph{The operator $L^{(r)}$ is then given by:
\begin{align*}
(L^{(r)}t)_{i_{1}\cdots i_{r}i_{r+1}i_{r+2}}^{k} &  =\frac{1}{2}\left(
t_{i_{1}\cdots i_{r}i_{r+1}i_{r+2}}^{k}-t_{i_{1}\cdots i_{r}i_{r+2}i_{r+1}%
}^{k}\right.  \\
&  +t_{i_{1}\cdots i_{r}ki_{r+1}}^{i_{r+2}}-t_{i_{1}\cdots i_{r}i_{r+1}%
k}^{i_{r+2}}\\
&  +\left.  t_{i_{1}\cdots i_{r}ki_{r+2}}^{i_{r+1}}-t_{i_{1}\cdots
i_{r}i_{r+2}k}^{i_{r+1}}\right)  .
\end{align*}
}

\emph{Taking the cyclic sum of this expression, and using (\ref{W1}), we
obtain:
\begin{align*}
\underset{\left(  i_{1}\cdots i_{r+1}\right)  }{\frak{S}}(L^{(r)}%
t)_{i_{1}\cdots i_{r+1}i_{r+2}}^{k} &  =\frac{1}{2}\left(  (r+1)t_{i_{1}\cdots
i_{r}i_{r+1}i_{r+2}}^{k}+t_{i_{1}\cdots i_{r}i_{r+1}i_{r+2}}^{k}\right.  \\
&  -t_{i_{1}\cdots i_{r}i_{r+1}k}^{i_{r+2}}-(r+1)t_{i_{1}\cdots i_{r}i_{r+1}%
k}^{i_{r+2}}\\
&  +\underset{\left(  i_{1}\cdots i_{r+1}\right)  }{\frak{S}}\left(  \left.
t_{i_{1}\cdots i_{r}ki_{r+2}}^{i_{r+1}}-t_{i_{1}\cdots i_{r}i_{r+2}k}%
^{i_{r+1}}\right)  \right)
\end{align*}
}

\emph{Therefore, equations (\ref{W2}) can be written:
\begin{equation}
t_{i_{1}\cdots i_{r}i_{r+1}i_{r+2}}^{k}-t_{i_{1}\cdots i_{r}i_{r+1}k}%
^{i_{r+2}}=\frac{1}{r+2}\underset{\left(  i_{1}\cdots i_{r+1}\right)
}{\frak{S}}\left(  t_{i_{1}\cdots i_{r}i_{r+2}k}^{i_{r+1}}-t_{i_{1}\cdots
i_{r}ki_{r+2}}^{i_{r+1}}\right)  .\label{previous}%
\end{equation}
}

\emph{Taking the cyclic sum with respect to $i_{1},\ldots,i_{r+2}$ on both
sides of last equation (using again (\ref{W1})) yields:
\begin{align*}
-\underset{\left(  i_{1}\cdots i_{r+2}\right)  }{\frak{S}}t_{i_{1}\cdots
i_{r}i_{r+1}k}^{i_{r+2}}=\frac{1}{r+2}\left(  (r+1)\underset{\left(
i_{1}\cdots i_{r+2}\right)  }{\frak{S}}t_{i_{1}\cdots i_{r}i_{r+1}k}^{i_{r+2}%
}+\underset{\left(  i_{1}\cdots i_{r+2}\right)  }{\frak{S}}t_{i_{1}\cdots
i_{r}i_{r+1}k}^{i_{r+2}}\right)  \\
=\underset{\left(  i_{1}\cdots i_{r+2}\right)  }{\frak{S}}t_{i_{1}\cdots
i_{r}i_{r+1}k}^{i_{r+2}},
\end{align*}
from which it follows that:
\[
\underset{\left(  i_{1}\cdots i_{r+2}\right)  }{\frak{S}}t_{i_{1}\cdots
i_{r}i_{r+1}k}^{i_{r+2}}=0.
\]
Last equation can be introduced in (\ref{previous}) giving:
\[
t_{i_{1}\cdots i_{r}i_{r+1}i_{r+2}}^{k}-t_{i_{1}\cdots i_{r}i_{r+1}k}%
^{i_{r+2}}=\frac{1}{r+2}\left(  -t_{i_{1}\cdots i_{r}i_{r+1}k}^{i_{r+2}%
}+t_{i_{1}\cdots i_{r}i_{r+1}i_{r+2}}^{k}\right)  ,
\]
from which one can conclude that the equations characterizing the subspace
$W^{r+1}\subset V_{1}^{r+1,1}$ are (\ref{W1}) and
\begin{equation}
t_{i_{1}\cdots i_{r}i_{r+1}i_{r+2}}^{k}-t_{i_{1}\cdots i_{r}i_{r+1}k}%
^{i_{r+2}}=0,\label{symmetric}%
\end{equation}
for any indices $1\leq k,i_{1},\ldots,i_{r+2}\leq n$. }

\emph{It should be pointed that, in \cite{E}, Epstein describes the space of
$\infty$-jets of Riemannian metrics $g$ at a point, as follows. He writes the
Taylor series for $g_{ij}(x)$ in normal coordinates: $\delta_{ij}+\sum
_{r\geq1}g_{iji_{1}\cdots i_{r}}x^{i_{1}}\cdots x^{i_{r}}$. The coefficients
$g_{iji_{1}\cdots i_{r}}$ satisfy the following conditions: }

\begin{enumerate}
\item \emph{\label{Sym1} They are symmetric in the first two indices. }

\item \emph{\label{Sym2} They are symmetric in the last $r$-indices. }

\item \emph{\label{Sym3} $\underset{\left(  i_{1}\cdots i_{r+1}\right)
}{\frak{S}}g_{ii_{1}\cdots i_{r+1}}=0$ for $1\leq i,i_{1},\ldots,i_{r+1}\leq
n$. }
\end{enumerate}

\emph{Then, he defines $f_{r}\in\left(  V^{\star}\right)  ^{\otimes(r+2)}$ by
$f_{r}(v_{i},v_{j},v_{i_{1}},\ldots,v_{i_{r}})=g_{iji_{1}\cdots i_{r}}$, and
proves that the set of $f_{r}\in\left(  V^{\star}\right)  ^{\otimes(r+2)}$
satisfying the symmetry conditions corresponding to \ref{Sym1}-\ref{Sym3} is
an irreducible $GL(n,\mathbb{R})$-module $Y_{r}$ with Young diagram having $r$
squares in the first row and $2$ squares in the second one, except that if
$r=1$ then $Y_{r}=\{0\}$. Moreover, for any sequence $f_{r}\in Y_{r}$, $2\leq
r<\infty$, there is a Riemannian metric whose Taylor series gives the elements
$f_{r}$, so that one can regard $\Pi_{r\geq2}Y_{r}$ as the space of $\infty
$-jets of Riemannian metrics at a point. }

\emph{The identification between Riemannian metrics $g$ and sections $s$ of
the bundle $FV/O(n)\rightarrow V$ together with Lemma \ref{uniqueness} yield
an isomorphism of $O(n)$-modules: $\mathbf{S}^{r}\cong\Pi_{2\leq k\leq r}%
Y_{k}$. Moreover, it is straightforward to check that the mapping
$g_{jki_{1}\cdots i_{r+1}}\mapsto t_{i_{1}\cdots i_{r}i_{r+1}j}^{k}$ also
defines an isomorphism of $O(n)$-modules $W^{r+1}\cong Y_{r+1}$. }
\end{ejemplo}

\begin{ejemplo}
[$O(p)\times O(q)$-structures]\emph{From the expression of $\delta^{-1}\circ
P_{\mathrm{im}\delta}$ given in Example \ref{foliations}, it is immediate that
the components of $L^{(r)}t$ are:
\begin{align}
(L^{(r)}t)_{i_{1}\cdots i_{r}i_{r+1}i_{r+2}}^{k}=\frac{1}{2}\left(
t_{i_{1}\cdots i_{r}i_{r+1}i_{r+2}}^{k}-t_{i_{1}\cdots i_{r}i_{r+2}i_{r+1}%
}^{k}\right.  \label{foliationsL1}\\
+t_{i_{1}\cdots i_{r}ki_{r+1}}^{i_{r+2}}-t_{i_{1}\cdots i_{r}i_{r+1}%
k}^{i_{r+2}}\nonumber\\
+\left.  t_{i_{1}\cdots i_{r}ki_{r+2}}^{i_{r+1}}-t_{i_{1}\cdots i_{r}i_{r+2}%
k}^{i_{r+1}}\right)  \nonumber
\end{align}
if $i_{r+1},i_{r+2},k\in I_{\alpha}$, $\alpha=1,2$,
\begin{align}
(L^{(r)}t)_{i_{1}\cdots i_{r}i_{r+1}i_{r+2}}^{k}=\frac{1}{2}\left(
t_{i_{1}\cdots i_{r}i_{r+1}i_{r+2}}^{k}-t_{i_{1}\cdots i_{r}i_{r+2}i_{r+1}%
}^{k}\right.  \label{foliationsL2}\\
-\left.  t_{i_{1}\cdots i_{r}i_{r+1}k}^{i_{r+2}}+t_{i_{1}\cdots i_{r}ki_{r+1}%
}^{i_{r+2}}\right)  \nonumber
\end{align}
if $i_{r+1}\in I_{\alpha}$, $i_{r+2},k\in I_{\beta}$, $\alpha\neq\beta$, and
\[
(L^{(r)}t)_{i_{1}\cdots i_{r}i_{r+1}i_{r+2}}^{k}=0
\]
if $i_{r+2}\in I_{\alpha}$, $k\in I_{\beta}$, $\alpha\neq\beta$. }

\emph{As in the previous example, relation (\ref{foliationsL1}) leads to the
equations (\ref{symmetric}) if $i_{r+1},i_{r+2},k\in I_{\alpha}$, $\alpha
=1,2$. On the other hand, taking the cyclic sum with respect to $i_{1}%
,\ldots,i_{r+1}$ in (\ref{foliationsL2}), and using (\ref{W1}), we obtain that
equations (\ref{symmetric}) are also satisfied if $i_{r+1}\in I_{\alpha}$,
$i_{r+2},k\in I_{\beta}$, $\alpha\neq\beta$. In this way, the equations of
$W^{r+1}\subset V_{1}^{r+1,1}$ turn out to be (\ref{W1}) and:
\begin{equation}
t_{i_{1}\cdots i_{r}i_{r+1}i_{r+2}}^{k}-t_{i_{1}\cdots i_{r}i_{r+1}k}%
^{i_{r+2}}=0,
\end{equation}
for any indices $i_{1},\ldots,i_{r+1}\in I_{1}\cup I_{2}$, $i_{r+2},k\in
I_{\alpha}$, $\alpha=1,2$. }
\end{ejemplo}

\begin{ejemplo}
[$\mathbb{R}^{\star}$-structures]\emph{The expression of $\delta^{-1}\circ
P_{\mathrm{im}\delta}$ given in Example \ref{webs} leads to the following
expression of $L^{(r)}$:
\[
(L^{(r)}t)_{i_{1}\cdots i_{r+1}i_{r+2}}^{k}=\frac{1}{n-1}\delta_{i_{r+2}}%
^{k}\sum_{j}\left(  t_{i_{1}\cdots i_{r+1}j}^{j}-t_{i_{1}\cdots ji_{r+1}}%
^{j}\right)  .
\]
}

\emph{Taking the cyclic sum with respect to $i_{1},\ldots,i_{r+1}$ (for
$k=i_{r+2}$), and using (\ref{W1}), we obtain:
\begin{align*}
\underset{\left(  i_{1}\cdots i_{r+1}\right)  }{\frak{S}}(L^{(r)}%
t)_{i_{1}\cdots i_{r+1}k}^{k}=\\
=\frac{1}{n-1}\left(  (r+1)\sum_{j}t_{i_{1}\cdots i_{r}i_{r+1}j}^{j}+\sum
_{j}t_{i_{1}\cdots i_{r}i_{r+1}j}^{j}\right)  =\\
=\frac{r+2}{n-1}\sum_{j}t_{i_{1}\cdots i_{r}i_{r+1}j}^{j}%
\end{align*}
}

\emph{The resulting equations of $W^{r+1}$ are then (\ref{W1}) and:
\begin{equation}
\sum_{k}t_{i_{1}\cdots i_{r+1}k}^{k}=0,\hspace{1cm}1\leq k,i_{1}%
,\ldots,i_{r+1}\leq n.
\end{equation}
}
\end{ejemplo}

\end{document}